\documentclass[preprint]{elsarticle}
\usepackage[colorlinks=true,linkcolor=blue]{hyperref}%

\usepackage{url}
\usepackage{amssymb, amsmath, amsthm} 
\usepackage{epsfig,graphicx} 
\usepackage{caption}      	%customized figure captions
\usepackage{subcaption}   	%for figures
\usepackage{color}
\usepackage{color}
\newtheorem*{defn}{Definition}

\newcommand{\bR}{{\mathbb{ R}}}

\newcommand{\cK}{{\cal K}}

\newcommand{\cW}{{\cal W}}

\newcommand{\eps}{\varepsilon}
\newcommand{\diam}[1]{\mathrm{diam}(#1)}

\newcommand{\Eq}[1]{(\ref{eqn:#1})}
\newcommand{\Fig}[1]{Fig.~\ref{fig:#1}}
\newcommand{\Sec}[1]{\S\ref{sec:#1}}

\newcommand{\beq}[1]{\begin{equation}\label{eqn:#1}}
\newcommand{\eeq}{\end{equation}}

\usepackage{amssymb}
\usepackage[figuresright]{rotating}

\begin{document}
\parskip=3pt

\begin{frontmatter}

%% Title, authors and addresses

%% use the tnoteref command within \title for footnotes;
%% use the tnotetext command for the associated footnote;
%% use the fnref command within \author or \address for footnotes;
%% use the fntext command for the associated footnote;
%% use the corref command within \author for corresponding author footnotes;
%% use the cortext command for the associated footnote;
%% use the ead command for the email address,
%% and the form \ead[url] for the home page:
%%
%% \title{Title\tnoteref{label1}}
%% \tnotetext[label1]{}
%% \author{Name\corref{cor1}\fnref{label2}}
%% \ead{email address}
%% \ead[url]{home page}
%% \fntext[label2]{}
%% \cortext[cor1]{}
%% \address{Address\fnref{label3}}
%% \fntext[label3]{}

%\dochead{}
%% Use \dochead if there is an article header, e.g. \dochead{Short communication}

\title{Exploring the Topology of Dynamical Reconstructions}

% \title{Exploring the \alert{(birth of /premature)} Topology of
%   \alert{(early/projected/compressed/under-embedded/improper)}Dynamical
%   Reconstructions}

%% use optional labels to link authors explicitly to addresses:
%% \author[label1,label2]{<author name>}
%% \address[label1]{<address>}
%% \address[label2]{<address>}

\author[cucs]{Joshua Garland}
\author[cucs,sfi]{Elizabeth Bradley}
\author[cuappm]{James D. Meiss}

\address[cucs]{Department of Computer Science, University of Colorado, Boulder CO USA}
\address[sfi]{Santa Fe Institute, Santa Fe NM USA}
\address[cuappm]{Department of Applied Mathematics, University of Colorado, Boulder CO USA}

\begin{abstract}

Computing the state-space topology of a dynamical system from scalar
data requires accurate reconstruction of those dynamics and
construction of an appropriate simplicial complex from the results.
The reconstruction process involves a number of free parameters and
the computation of homology for a large number of simplices can be
expensive.  This paper is a study of how to avoid a full
(diffeomorphic) reconstruction and how to decrease the computational
burden.  Using trajectories from the classic Lorenz system, we
reconstruct the dynamics using the method of delays, then build a
parsimonious simplicial complex---the ``witness complex''---to compute
its homology.  Surprisingly, we find that the witness complex
correctly resolves the homology of the underlying invariant set from
noisy samples of that set even if the reconstruction dimension is well
below the thresholds specified in the embedding theorems for assuring
topological conjugacy between the true and reconstructed dynamics.  We
conjecture that this is because the requirements for reconstructing
homology, are less stringent than those in the embedding theorems.  In
particular, to faithfully reconstruct the homology, a homeomorphism is
sufficient---as opposed to a diffeomorphism, as is necessary for the
full dynamics.  We provide preliminary evidence that a homeomorphism,
in the form of a delay-coordinate reconstruction map, may manifest at
a lower dimension than that required to achieve an embedding.

\end{abstract}

\begin{keyword}
%% keywords here, in the form: keyword \sep keyword
Topology \sep Delay-Coordinate Embedding\sep Nonlinear Time Series Analysis\sep Computational Homology \sep Witness Complex
%% MSC codes here, in the form: \MSC code \sep code
%% or \MSC[2008] code \sep code (2000 is the default)

\end{keyword}

\end{frontmatter}

%%
%% Start line numbering here if you want
%%
% \linenumbers%

%%%%%%%%%%%%%%%%%%%%%%%
\section{Introduction}\label{sec:intro}

Topology is of particular interest in dynamics, since many
properties---the existence of periodic orbits, transitivity,
recurrence, entropy, etc.---depend only upon topology. This idea is
commonly exploited in the computational topology community, often
using the Conley index of isolating neighborhoods, to study dynamical
invariants \cite{Kaczynski04}.  However, computing topology from time
series can be a real challenge.  First, one typically has only {\sl
  scalar} data, not the full trajectory, and hence one must begin by
reconstructing the full dynamics from that data---e.g., via
delay-coordinate reconstruction.  Success of this reconstruction
procedure depends on several free parameters.  In practice, the
embedding theorems provide little guidance regarding how to choose
these parameters; though a number of creative strategies have been
developed, these methods require good data and input from a human
expert.  Moreover, the delay-coordinate reconstruction machinery (both
theorems and heuristics) targets the computation of dynamical
invariants like the correlation dimension and the Lyapunov exponent.
If one just wants to extract the topological structure of an invariant
set, as we show in this paper, that full machinery may not be needed.
Nevertheless, there are scales and scale parameter issues here, as for
the standard machinery.  Moreover, real-world data sets have finite
length, sample some underlying set at a finite time interval, have
finite precision, and may be contaminated by noise.  In the face of
these issues, one obviously cannot compute the topology to arbitrary
precision.

Coarse-graining the topological analysis of data also addresses
another issue: the associated computations are expensive, and that
expense grows with the number of simplices in the complexes that one
constructs during that process.  The pioneering work in this area used
cubical complexes and multivalued maps for this purpose
\cite{mischaikow99}, and these results can be computationally rigorous
even in the face of noise.  For more efficiency, one can use a
simplicial complex that follows the natural geometry of the
data---e.g., the {\sl witness complex} of \cite{deSilva04}.  To
construct a witness complex, one chooses a set of ``landmarks,''
typically a subset of the data, that become the vertices of the
complex.  The connections between the landmarks are determined by their
nearness to the rest of the data---the ``witnesses.''  Two landmarks
in the complex are joined by an edge, for instance, if they share at
least one witness.  As described in \Sec{wc}, there are many possible
definitions for a witness ``relation."  The one that we use includes a
scale parameter, $\eps$, intended to provide a measure of noise
immunity.  The ideas of persistent homology \cite{ELZ01,Robins02} can
be used to choose $\eps$, build the complex, and then explore the
changes in its topology with changing reconstruction dimension.  

Our initial work on this approach suggests that, for a number of
different dynamical systems, \emph{the witness complex correctly
  resolves the homology of the underlying invariant set
% ---viz., its Betti numbers---
even if the reconstruction dimension is well below the thresholds for
which the embedding theorems assure smooth conjugacy between the true
and reconstructed dynamics.}  This paper reports upon an exploration
of that conjecture in the context of the classic Lorenz system and
suggests some implications and applications.  To set the stage for
that discussion, the rest of this section gives a brief review of
delay-coordinate reconstruction.
% 
% : that one can compute the topology of an attractor of a dynamical
% system using a coarse-grained simplicial complex built from a
% low-dimensional reconstruction of scalar data sampled from that
% attractor.
% 
The witness complex is covered in more depth in \Sec{wc}, which also
describes the notion of persistence and demonstrates how that idea is
used to choose scale parameters for a complex built from reconstructed
time-series data. In \Sec{meat}, we explore how the homology of such a
complex changes with reconstruction dimension.

% \subsection{Delay-coordinate reconstruction} 

\begin{figure}
\begin{centering}
        \begin{subfigure}[b]{0.32\textwidth}
                \includegraphics[width=\textwidth]{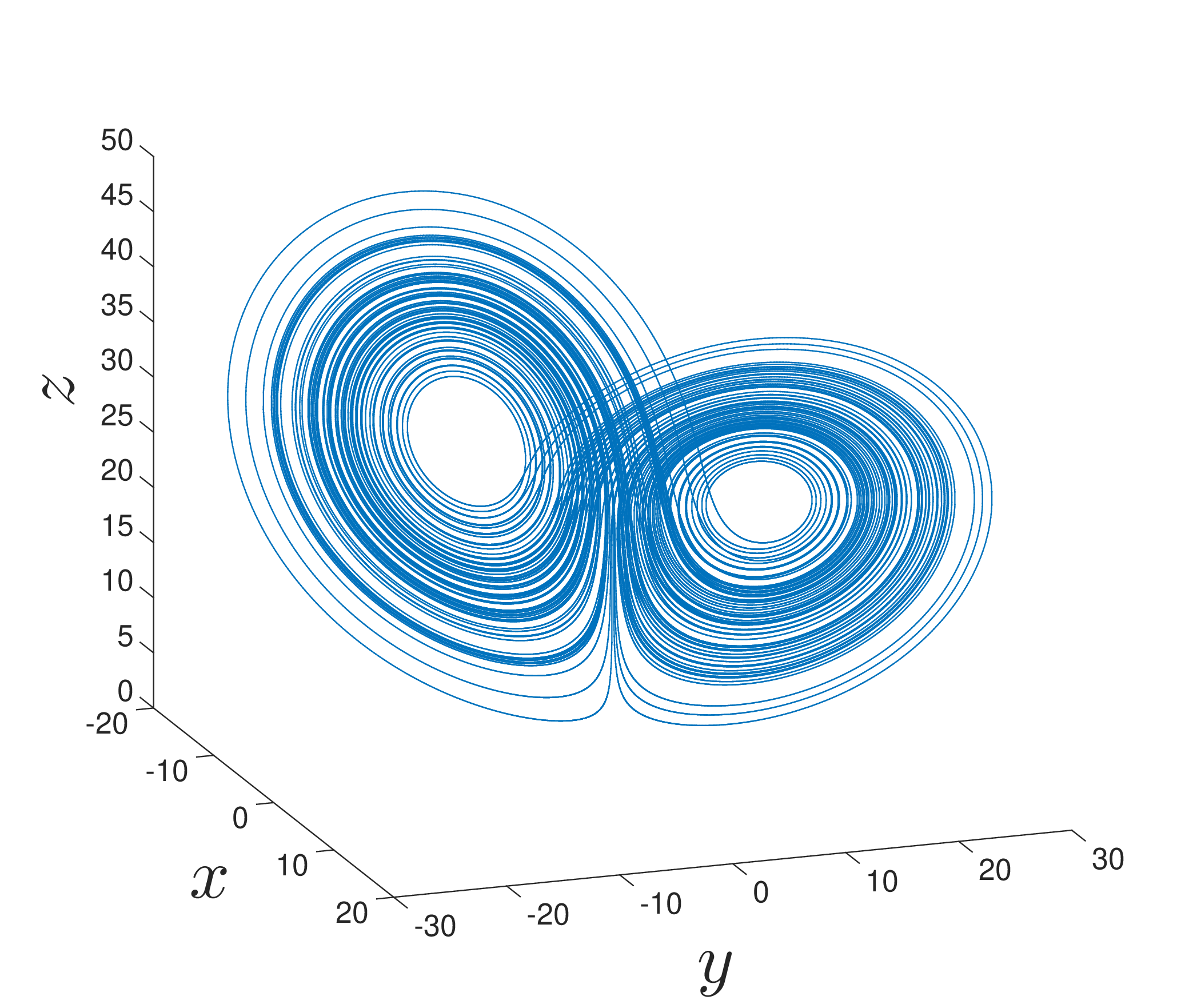}
                \caption{}
        \end{subfigure}
        \begin{subfigure}[b]{0.32\textwidth}
                \includegraphics[width=\textwidth]{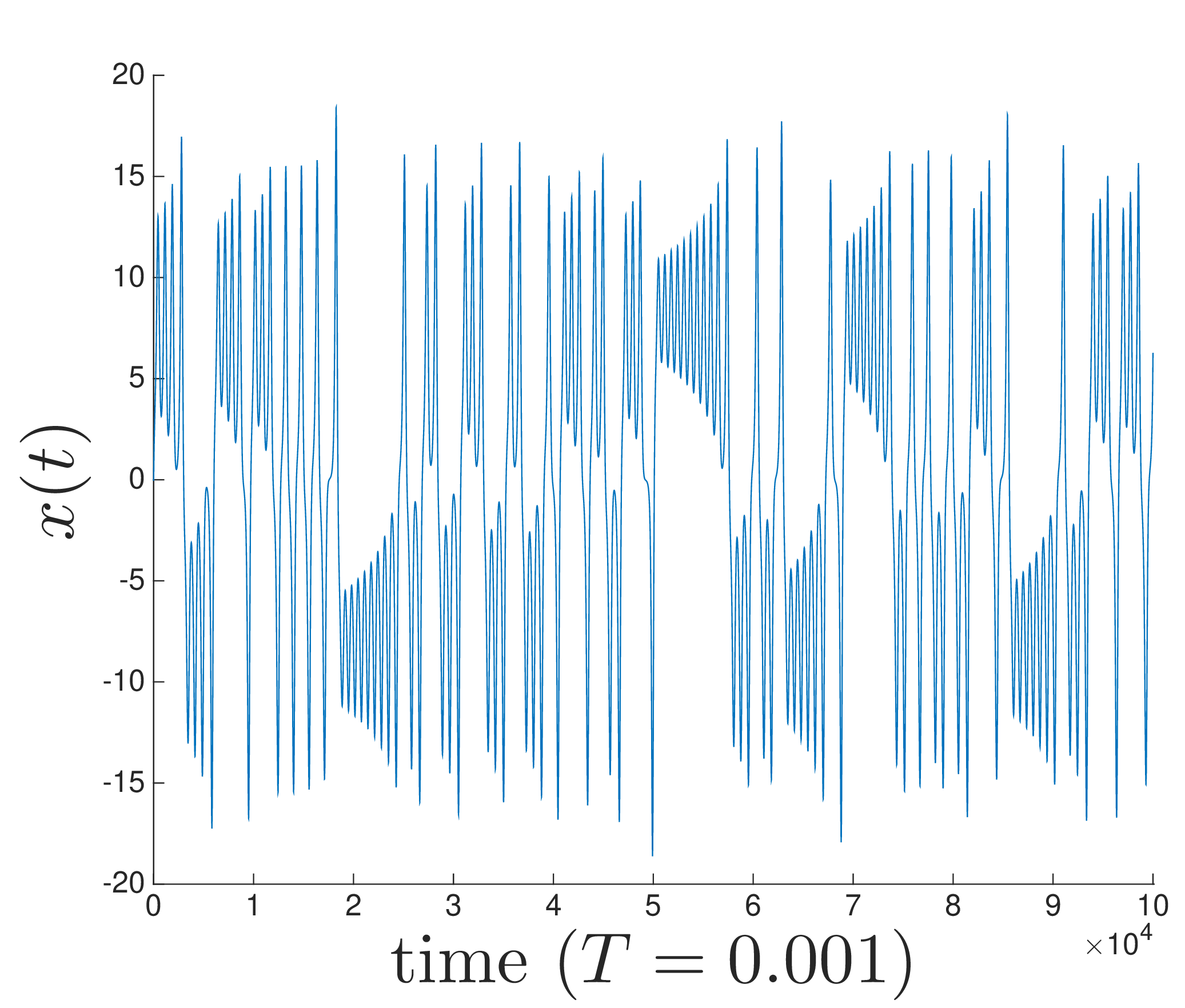}
                \caption{}
        \end{subfigure}
        \begin{subfigure}[b]{0.32\textwidth}
                \includegraphics[width=\textwidth]{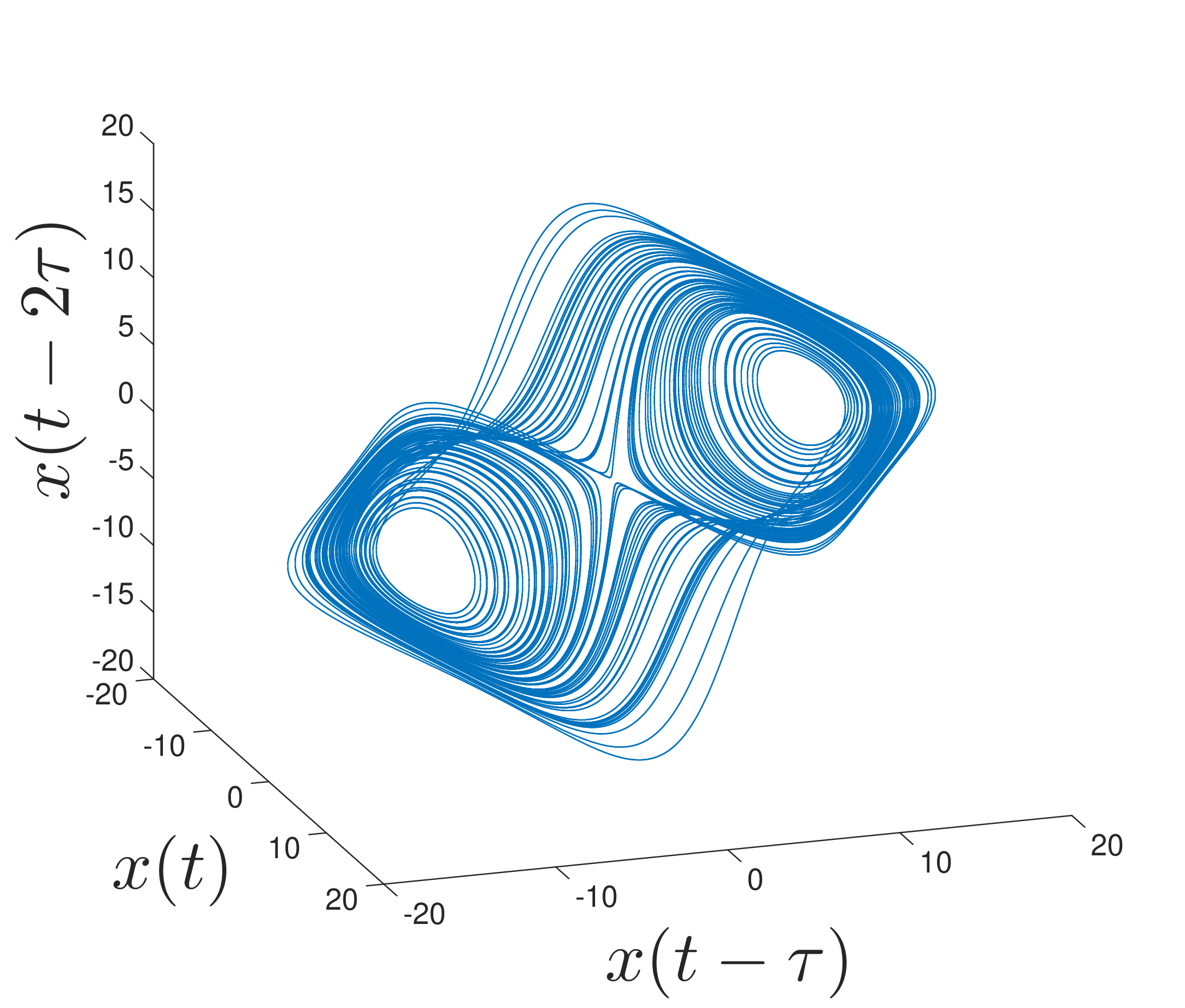}
                \caption{}
        \end{subfigure}
        \end{centering}
\caption{Classic Lorenz attractor ($r=28$, $b=8/3$, $\sigma=10$): (a)
  A $10^5$-point trajectory in $\bR^3$ generated using fourth-order
  Runga-Kutta with a time step of $T=0.001$.  (b) A time-series trace
  of the $x$ coordinate of that trajectory.  (c) A $3D$ projection of
  a delay-coordinate embedding with dimension $m=5$ and delay $\tau =
  174 T$, following \Eq{takens}.}
\label{fig:Lorenz63Example}
\end{figure}

Delay-coordinate reconstruction \cite{packard80,Takens81} is arguably the
most well-established technique for reconstructing the dynamics of a
system from scalar time-series data.  Suppose that $\vec{Y}$ is a point on
a compact invariant set $M \subset \bR^d$, and $\vec{Y}(t)$ represents its trajectory.
A smooth {\sl measurement function}
$h:M \to \bR$ gives rise to a scalar time-series, $x(t) = h(\vec{Y}(t))$, from
that trajectory. Then the delay-coordinate map, $F: M \to \bR^m$
 %\vspace*{-1mm}
 % 
\begin{equation}\label{eqn:takens}
 F(\vec{Y}(t);h,m,\tau) = (x(t), ~ x(t-\tau), ~ \dots , ~x(t-(m-1)\tau)) \enspace,
 \end{equation}
is almost always a diffeomorphism whenever $\tau>0$ and if
$m$ is large enough, i.e.,
% $m > 2d$, according to Takens \cite{Takens81}, or 
$m>2d_{box}$, 
%according to Sauer {\sl et al.}, 
where $d_{box}$ is the box-counting dimension of $M$ \cite{sauer91}.
When these conditions are met, the reconstructed attractor and the
true attractor are diffeomorphic, and thus certainly have the same
topology.  The left panel of \Fig{Lorenz63Example} shows an example: a
trajectory from the classic Lorenz system \cite{Lorenz63}.  The middle
panel shows the corresponding time series of the $x$ coordinate of
that trajectory (i.e., $h(x,y,z) = x$), and the right panel shows a
delay-coordinate reconstruction using $\tau=174T$, where $T$ is the
interval between points in the time series.  Note that a
reconstruction dimension of five ($m=5$) is required in order to
satisfy the $m>2d_{box}$ requirement for this attractor, since
$d_{box} \approx 2.06$.  Of course, it is not easy to display the $5D$
picture; \Fig{Lorenz63Example}(c) shows a $3D$ projection of this
reconstruction.

%%%%%%%%%%%%%%%%
%%%%%%%%%%%%%%%%%%

In practice, one is presented with a scalar time series so that the
dimension $d$ of the original state space is unknown, and one cannot
compute $d_{box}$ without first embedding the data.  Thus, choosing
the reconstruction dimension $m$ is a challenge.  There are a number
of heuristics for doing so.  Perhaps the most well-known is the family
of false near-neighbor methods pioneered in \cite{KBA92}.  The basic
idea behind this class of methods is to increase the reconstruction
dimension until the geometry of the neighbor relationships stabilizes;
this is taken to indicate that any false crossings created by the
measurement function $h$ have been eliminated and the dynamics are
properly unfolded.  The choice of the delay $\tau$ also plays a role
in this unfolding.  Though the theorems only require $\tau>0$, in
practice one needs to ensure that $\tau$ is large enough to make the
coordinates numerically independent, but not so large that the
coordinates become causally unrelated \cite{Casdagli:1991a}.  The
standard approach for this---which we used to select the $\tau$ value
in \Fig{Lorenz63Example}(c)---is to calculate the time-delayed average
mutual information of the time series and choose $\tau$ at the first
minimum of that curve \cite{fraser-swinney}.  There are many other
methods for estimating both $m$ and $\tau$; see \cite{Bradley-Kantz15}
for a deeper discussion.  All of these procedures are subtle and
subjective.  Invoking them and interpreting their results
requires good data and expert knowledge; the false-near neighbor
method, for instance, has been shown to regularly overestimate
embedding dimension when noise is present in the time
series---something that is unavoidable in experimental data.

In this paper, we adopt the philosophy that one might only desire
knowledge of the topology of the invariant set, and we conjecture that
this might be possible with a lower reconstruction dimension than that
needed to obtain a true ``embedding.''  That is, the reconstructed
dynamics might be \emph{homeomorphic} to the original dynamics at a
lower dimension than that needed for a diffeomorphically correct
embedding.  We will return to this idea below.
%in \Sec{concl}.

%%%%%%%%%%%%%%%%%%%%%%%
\section{Witness Complexes for Dynamical Systems}
\label{sec:wc}

To compute the topology of a data set that samples an invariant set of
a dynamical system, we need a complex that captures the shape of the
data, but is robust with respect to noise and other sampling issues.
To do so \emph{efficiently}, the complex should be parsimonious.  A
witness complex is an ideal choice for these purposes.  Such a complex
is determined by the reconstructed time-series data, $W \subset
\bR^m$---the \emph{witnesses}---and an associated set $L \subset
\bR^m$, the \emph{landmarks}, which can (but need not) be chosen from
among the witnesses.
% 
% \footnote{Note that in order to do the false-witness stuff, one must
%   use the same lifting operation on all of the points involved in the
%   operation.  If the landmark set is not a subset of the data, this is
%   problematic.}.
% 
The landmarks form the vertex set of the complex; the connections
between them are dictated by the geometric relationships between $W$
and $L$.  In a general sense, a witness complex can be defined through
a relation $R(W,L) \subset W \times L$. As Dowker noted
\cite{Dowker52}, any relation gives rise to a pair of simplicial
complexes.  We will use one: a point $w \in W$ is a witness to an
abstract $k$-dimensional simplex $\sigma = \{l_{i_1}, l_{i_2}, \ldots
l_{i_{k+1}}\} \subset L$ whenever $\{w\} \times \sigma \subset
R(W,L)$. The collection of simplices that have witnesses is a complex
relative to the relation $R$. For example, two landmarks are connected
if they have a common witness---this is a one-simplex.  Similarly, if
three landmarks have a common witness, they form a two-simplex, and so
on.

There are many possible definitions for a witness relation $R$.  One
very natural construction is to use the matrix $D(W,L)$ of distances
$D_{ij} = \|w_i - l_j\|$ to define $R$. Sorting each row of this
matrix from smallest to largest determines the set of landmarks
that are closest to the $i^{th}$ witness. A relation corresponds to assigning a
cut-off, which thereby determines the simplices witnessed by
$w_i$. For example one can choose a fixed number (viz., $k$-nearest
neighbors), a strict size (neighbors within some distance), or an
increment.  The first concept gives the ``weak witness complex" of de
Silva and Carlsson \cite{deSilva04}, but suffers from the problem that
there is no limit on the distance to the nearest neighbors and thus a
simplex might be too spread out. The second notion seems too
restrictive: a portion of the invariant set $M$ that has a low density
may not be covered enough to be represented in the complex.  The third
idea is a compromise and gives the notion of an $\eps$-weak witness \cite{Carlsson14},
or what we call a ``fuzzy" witness \cite{Alexander15}: a point witnesses a simplex 
if all the landmarks in that simplex are within $\eps$ of the closest landmark to
the witness:
%%%%%
\begin{defn}[Fuzzy Witness]
The fuzzy witness set for a point $l \in L$ is the set
of witnesses
\beq{witnessSet}
	\cW_\eps(l) = \{w \in W : \|w-l\| \le \min_{l' \in L} \|w-l'\| + \eps \}.
\eeq
\end{defn}
%%%%%
\noindent In this case the relation consists of the collections
$R = \cup_{l \in L} (\cW_\eps(l) \times \{ l \})$, and a simplex $\sigma$ is
in the complex whenever $\cap_{l \in \sigma} \cW_\eps(l) \neq \emptyset$, that
is, all of its vertices share a witness.

The fuzzy witness complex reduces to ``strong witness complex'' of de
Silva and Carlsson \cite{deSilva04} when $\eps =0$.  In such a
complex, an edge exists between two landmarks \emph{iff} there exists
a witness that is exactly equidistant from those landmarks.  This is
not a practical notion of shared closeness for finite data
sets.  A simpler implementation of the fuzzy witness complex gives a
``clique" or ``flag" complex, analogous to the Rips complex
\cite{Ghrist08}, that consists of simplices whose pairs of vertices
have a common witness (this is called a ``lazy" complex in
\cite{deSilva04}). In this case the complex is
\beq{fuzzyComplex}
	\cK_\eps(W,L) = \{ \sigma \subset L: \cW_\eps(l) \cap \cW_\eps(l') \neq \emptyset, \;\; \forall l,l' \in \sigma\}
\eeq
Following \cite{Alexander15}, we will use this particular construction
because it minimizes computational complexity.

\begin{figure}[t]
        \begin{subfigure}[b]{0.45\textwidth}
                \includegraphics[width=\textwidth]{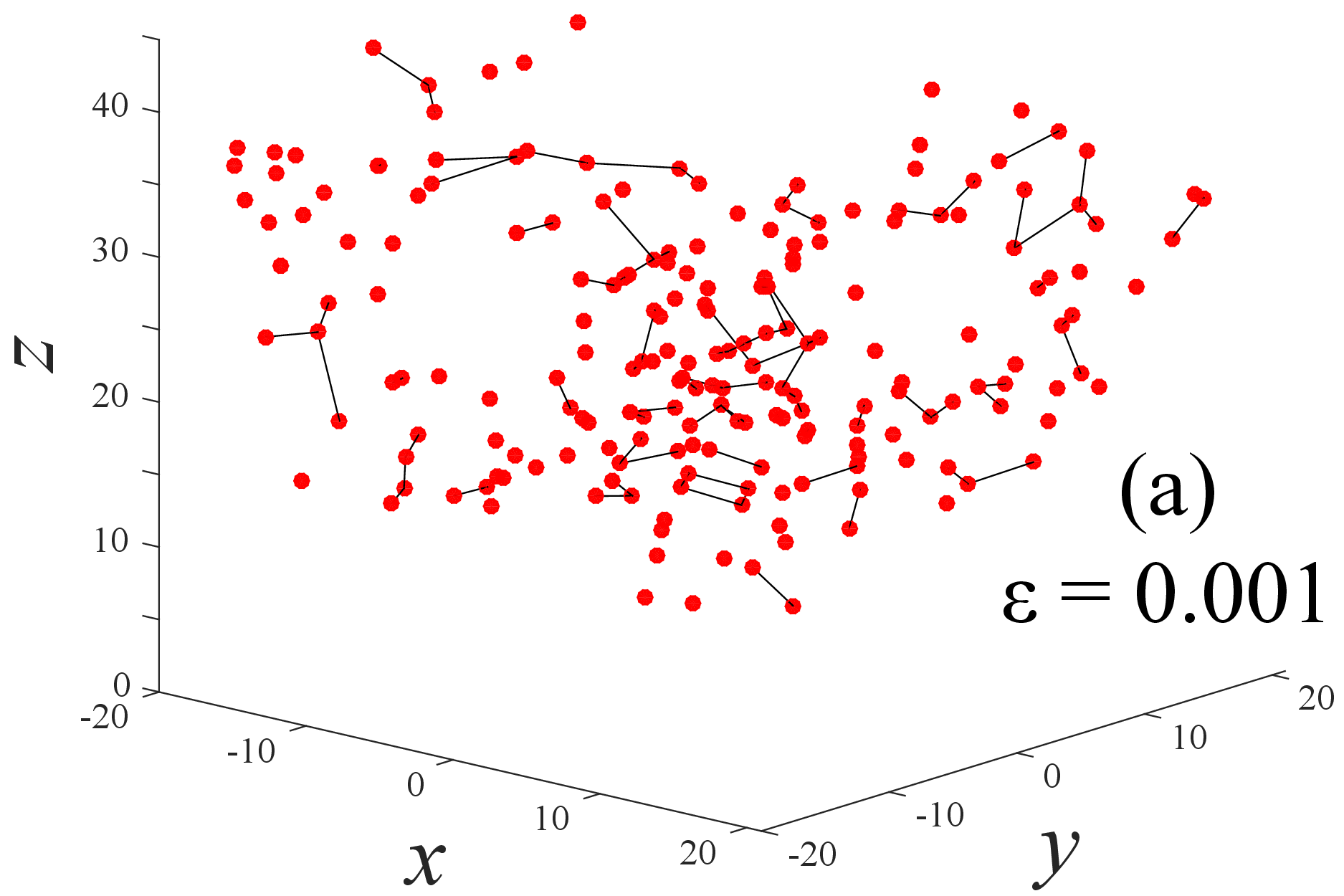}
%                \caption{$\eps=0.001$}
        \end{subfigure}
        \begin{subfigure}[b]{0.45\textwidth}
                \includegraphics[width=\textwidth]{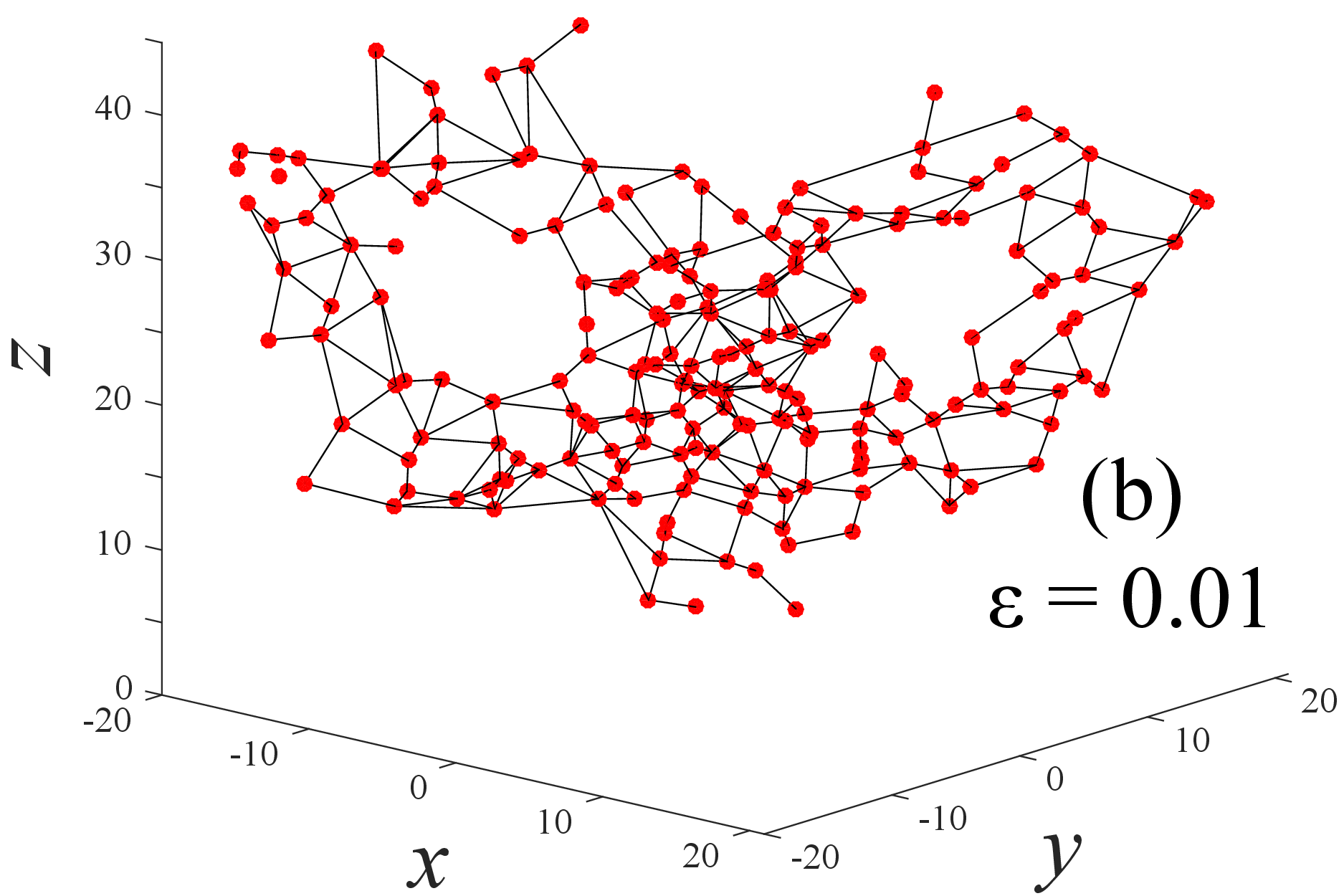}
%                \caption{$\eps=0.01$}
        \end{subfigure}
        \begin{subfigure}[b]{0.45\textwidth}
                \includegraphics[width=\textwidth]{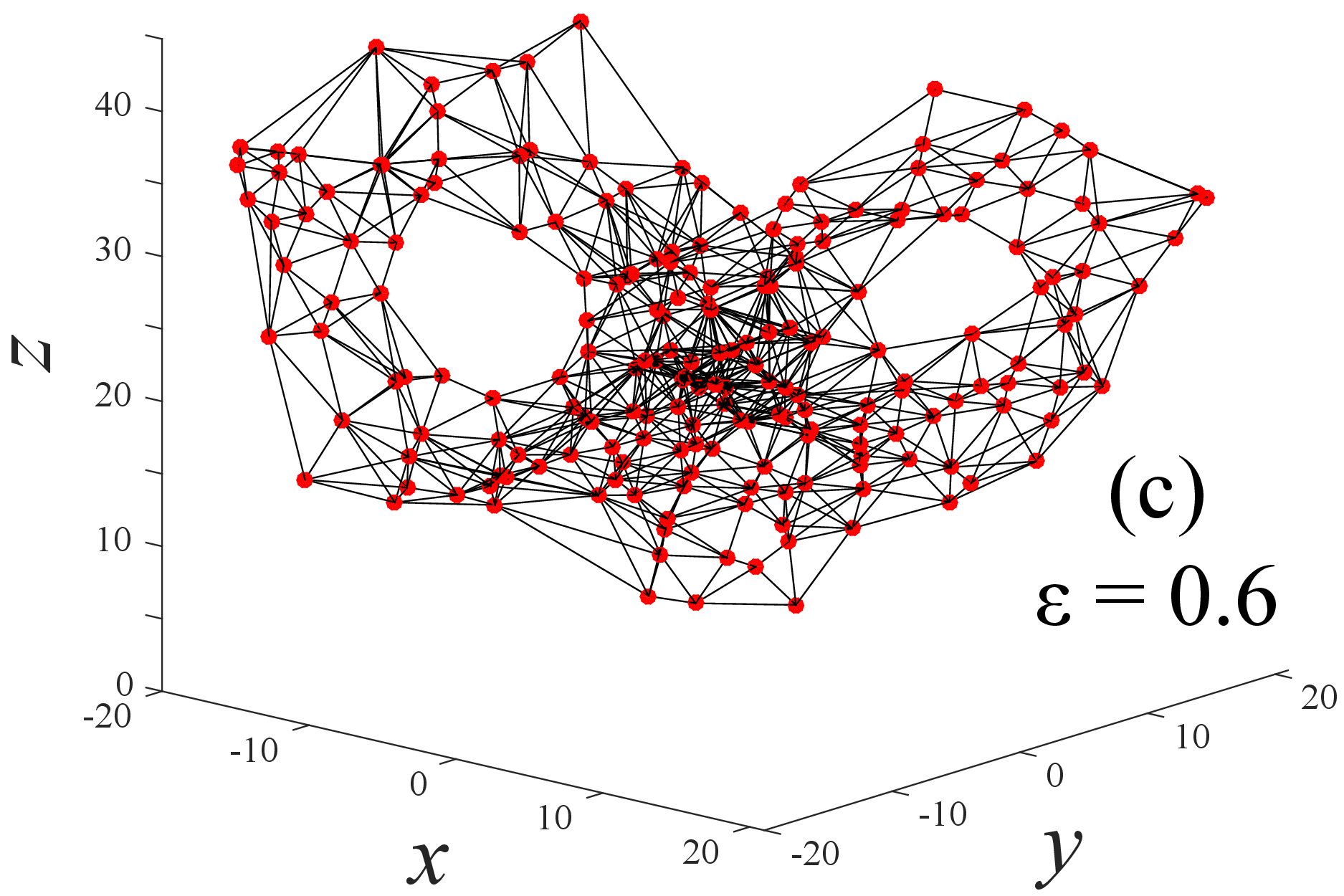}
%                \caption{$\eps=0.6$}
        \end{subfigure} \quad \quad \quad \quad
        \begin{subfigure}[b]{0.45\textwidth}
                \includegraphics[width=\textwidth]{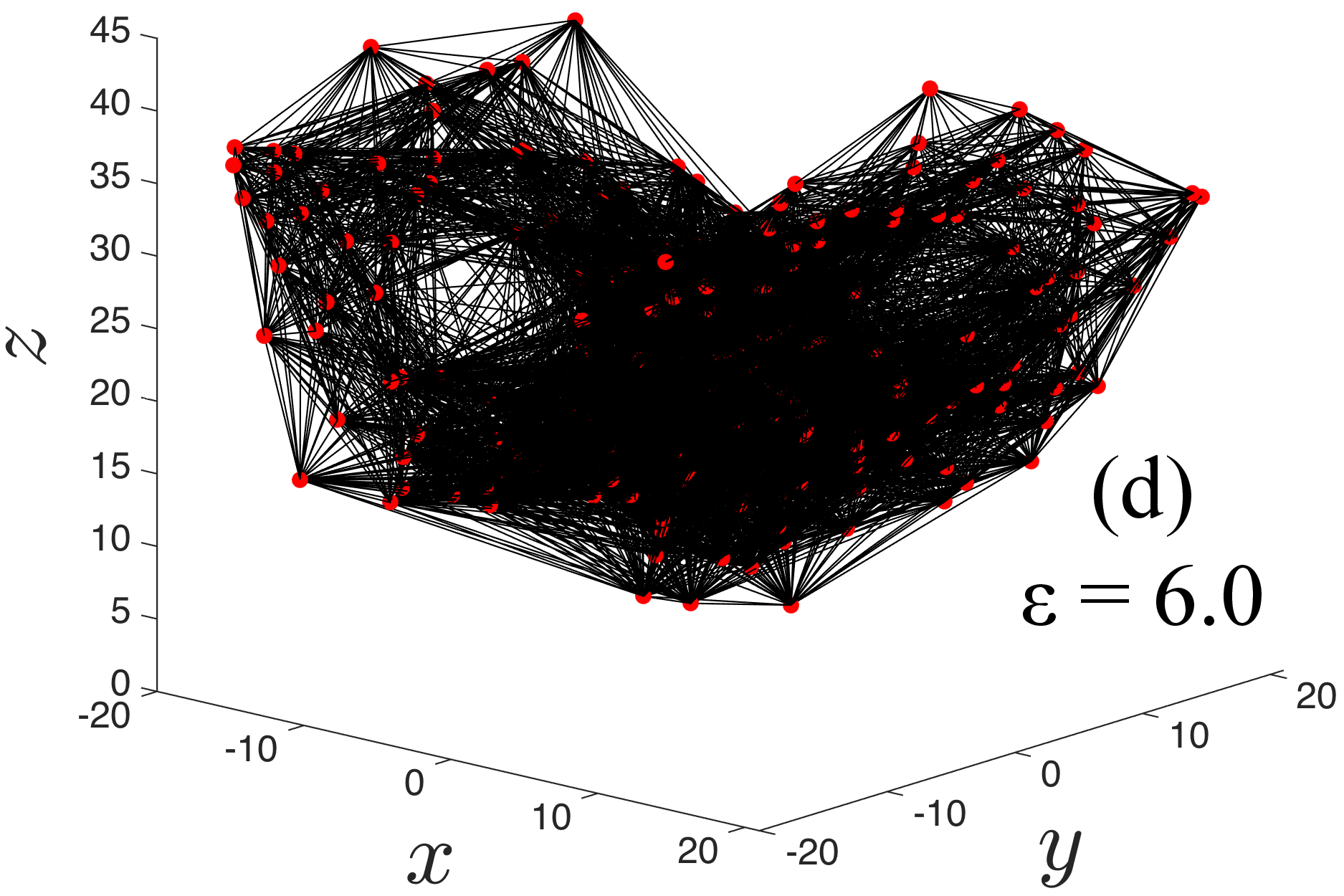}
%                \caption{$\eps=6$}
        \end{subfigure}
\caption{Varying the fuzziness parameter $\eps$: One-skeletons of
   clique complexes $\cK_\eps(W,L)$ constructed from the trajectory of
  \Fig{Lorenz63Example}(a) using 201 landmarks (red dots) and four
  values of $\eps$.
% for our reference:
%(a)$xi_\eps = 1.3286x10^{-5} = \eps/dA = 0.001/75.2668$ (b)
%$\xin = 1.3286x10^{-4}$ (c) $\xi = 0.0080$ (d)
%$\xi =0.0797$
}
\label{fig:witness-complexes-different-eps}
\end{figure}

Figure~\ref{fig:witness-complexes-different-eps} shows four witness
complexes built in this fashion from the 100,000-point trajectory of
the Lorenz system that is shown in \Fig{Lorenz63Example}(a).  The landmarks
(red dots) consist of $\ell = 201$ points equally spaced along the
trajectory, i.e., every $\Delta t = 500^{th}$ point of the time
series.  When $\eps$ is small, very few witnesses fall in the thin
regions required by \Eq{witnessSet}, so the resulting complex does not have
many edges and is thus not a good representation of the shape of the data. As $\eps$
grows, more witnesses fall in the ``shared'' regions and the complex
fills in, revealing the basic homology of the attractor of which the
trajectory is a sample.  There is an obvious limit to this, however:
when $\eps$ is very large, the even the largest holes in the complex are obscured.

%%%%%%
%%%%%

\begin{figure}
\centering
        \begin{subfigure}[b]{0.45\textwidth}
                \includegraphics[width=\textwidth]{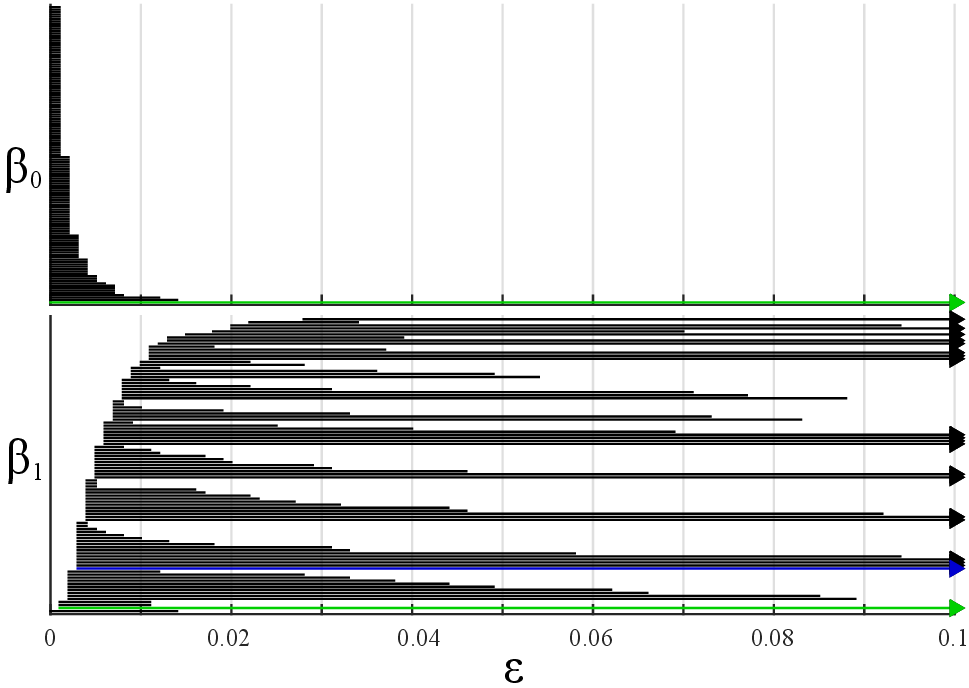}
        \end{subfigure}
        \begin{subfigure}[b]{0.45\textwidth}
			\includegraphics[width=\textwidth]{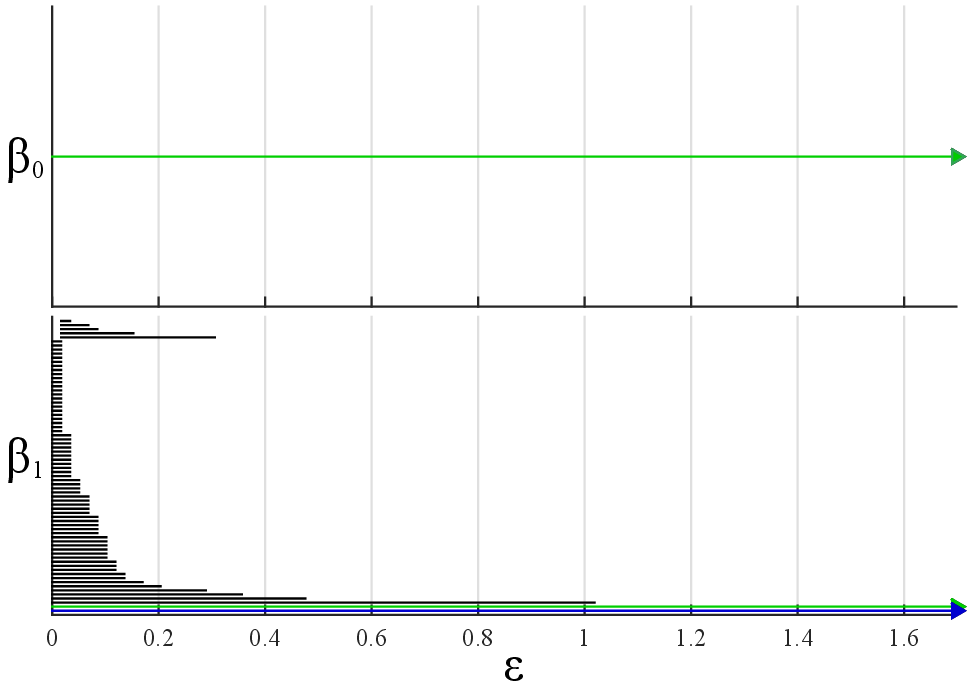}
        \end{subfigure}		
\caption{Persistence barcodes computed using {\tt javaPlex} for a
  $\ell=201$ witness complex of the trajectory of
  \Fig{Lorenz63Example}(a).  Each plot tabulates the two lowest Betti
  numbers of the complex for $100$ values of the scale parameter
  $\eps$.  The left panel shows the behavior when $0.001\le \eps\le
  0.1$, and the right $0.017 \le \eps \le 1.7$.}
 
\label{fig:Lorenz63-barcodes}
\end{figure}

Studying the change in homology under changing scale parameters is a
well-established notion in computational topology.  The underlying
idea of \emph{persistence} \cite{ELZ01,Robins02,Zomorodian05} is that
any topological property of physical interest should be (relatively)
independent of parameter choices in the associated algorithms.  One
useful way to represent information about the changing topology of a
complex is the \emph{barcode persistence diagram} \cite{Ghrist08}.
\Fig{Lorenz63-barcodes} shows barcodes of the first two
Betti numbers for the witness complexes of
\Fig{witness-complexes-different-eps}.  Each horizontal line in the
barcode is the interval in $\eps$ for which there exists a particular non-bounding
cycle, thus the number of such lines is the rank of the homology
group---a Betti number.  We computed these values for $\beta_0$ and
$\beta_1$ using {\tt javaPlex} \cite{javaplex} over the range $0.017\le \eps \le
1.7$, using the ``explicit" landmark selector to choose the equally
spaced points and the ``lazy" witness stream to obtain a clique
complex from the $\ell=201$ landmarks.  There were no three-dimensional voids, i.e., $\beta_2$ was 
always zero for this range of $\eps$---a reasonable implication for 
the $2.06$-dimensional attractor.  When $\eps$ is very small,
as in \Fig{witness-complexes-different-eps}(a), the witness complex
has many components and the $\beta_0$ barcode shows a large number of
entries. As $\eps$ grows, the spurious gaps between these components
disappear, leaving a single component that persists above $\eps
\approx 0.014$.  That is, witness complexes constructed with $\eps >
0.014$ correctly capture the connectedness of the underlying
attractor.  The $\beta_1$ barcode plot shows a similar pattern: there
are many holes for small $\eps$ that are successively filled in as
that parameter grows, leaving the two main holes (i.e., $\beta_1=2$)
for $\eps > 1.01$.  Above $\eps > 3.2$ (not shown in
\Fig{Lorenz63-barcodes}), one of those holes disappears; eventually,
for $\eps> 4.05$, the complex becomes topologically trivial.  Above
this value, the resulting complexes---recall
\Fig{witness-complexes-different-eps}(d)---have no non-contractible
loops and are homologous to a point (acyclic).

%%%%%
%%%%%%%

This notion of persistence can be turned around and used to select
good values for the parameters that play a role in topological data
analysis---i.e., looking for the $\eps$ value at which the homology
stabilizes.  However, definitions of what constitutes stabilization
are subjective and can be problematic.  This kind of issue turns up
routinely in nonlinear time-series analysis \cite{Bradley-Kantz15}.
Even so, persistence is a powerful technique and we make
use of it in a number of ways here.

Another critical step in our approach is the selection of the
landmarks that constitute the vertex set of the witness complex.
% 
% Since our goal is to explore changes in its topology with changing
% reconstruction dimension, we have to choose landmarks from among the data points.
% 
For efficiency, the number of landmarks should be much smaller than
the number of points in the time series, but for efficacy they should
be distributed so as to capture the shape of the data. A simple method
for this, advocated by \cite{deSilva04}, is to use a max-min algorithm
that chooses $L \subset W$ by selecting the farthest point in $W$ from
the previous selection and iterating until desired density and
sparseness requirements are satisfied, if possible. For data from a
dynamical system, one can alternatively exploit the natural temporal
ordering and select points that are equally spaced in time ($\Delta
t$) along the trajectory.  If the attractor is ergodic, this will
distribute the landmarks evenly relative its invariant measure.  One
advantage of this strategy in the context of delay-coordinate
reconstruction, as discussed in \Sec{meat} below, is that it allows
both witnesses and landmarks to be consistently moved from one
embedding dimension to another.  Thus we adopt this approach here.
% [[could add a footnote here that W is really W/L
Needless to say, the choice of the time interval between landmarks
invokes an tradeoff between accuracy and computational efficiency.

\begin{figure}
        \centering
        \begin{subfigure}[b]{0.3\textwidth}
                \includegraphics[width=\textwidth]{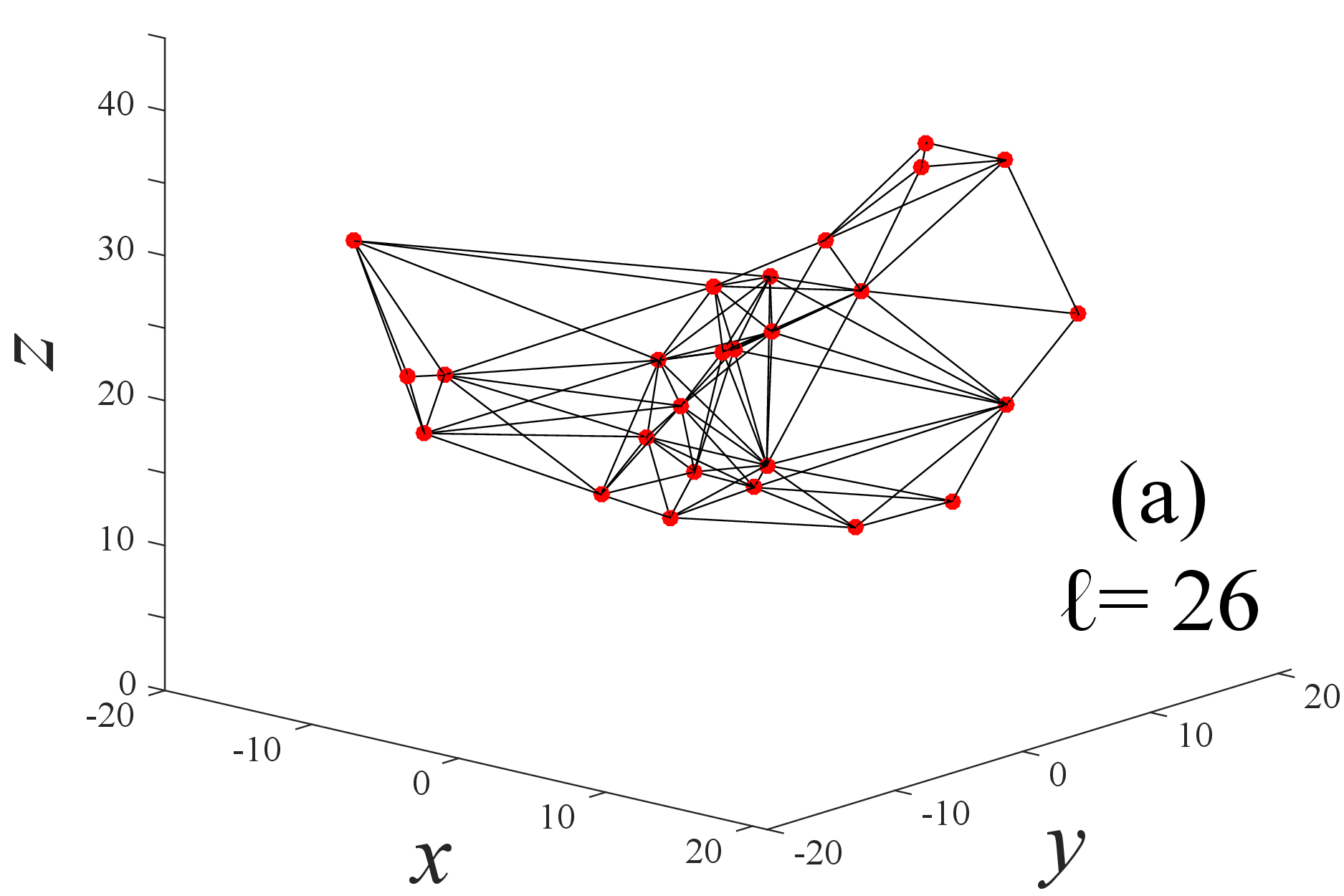}
  %               \caption{$l=26$}
        \end{subfigure}
         \begin{subfigure}[b]{0.3\textwidth}
                \includegraphics[width=\textwidth]{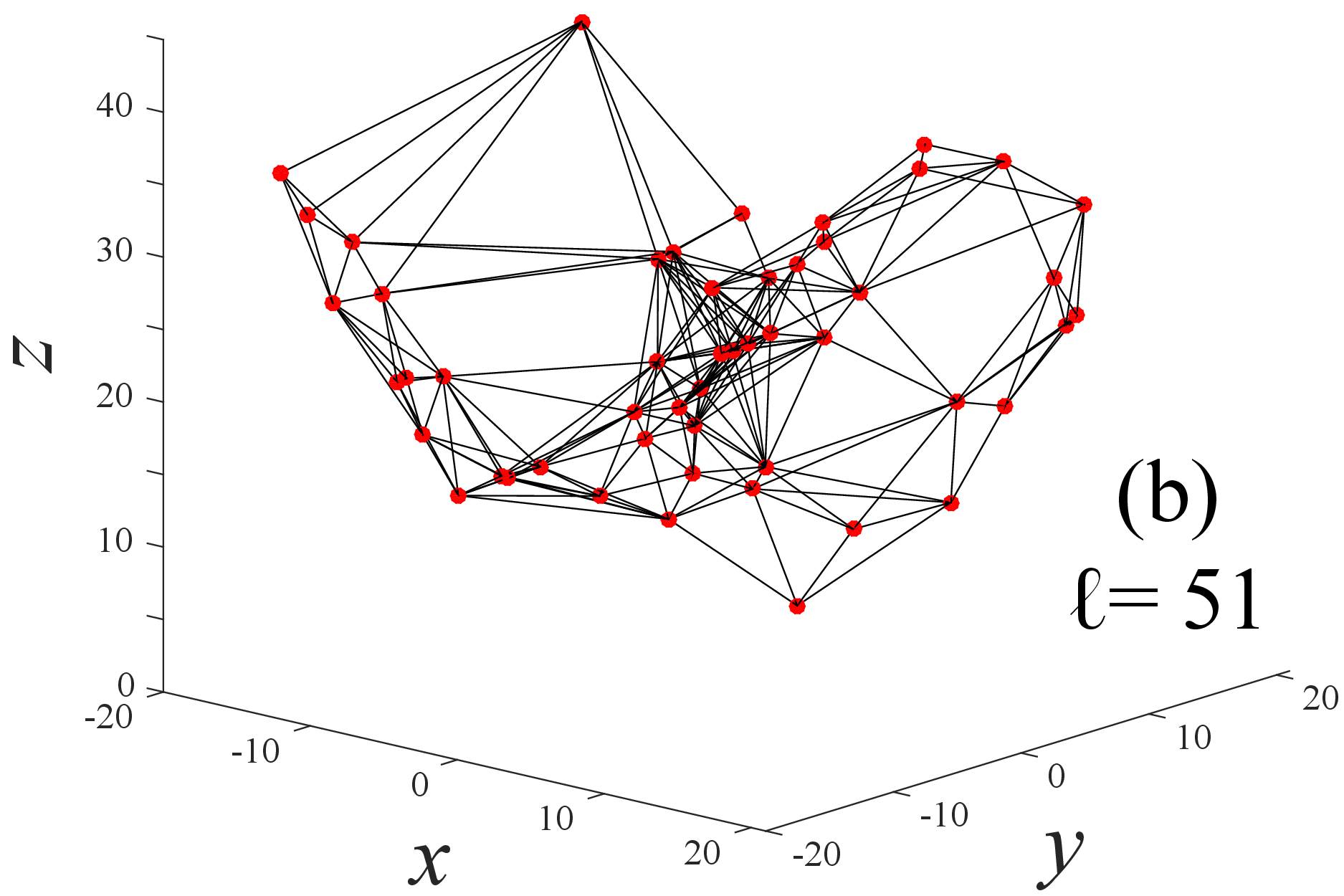}
   %              \caption{$l=51$}
                        \end{subfigure}
        \begin{subfigure}[b]{0.3\textwidth}
                \includegraphics[width=\textwidth]{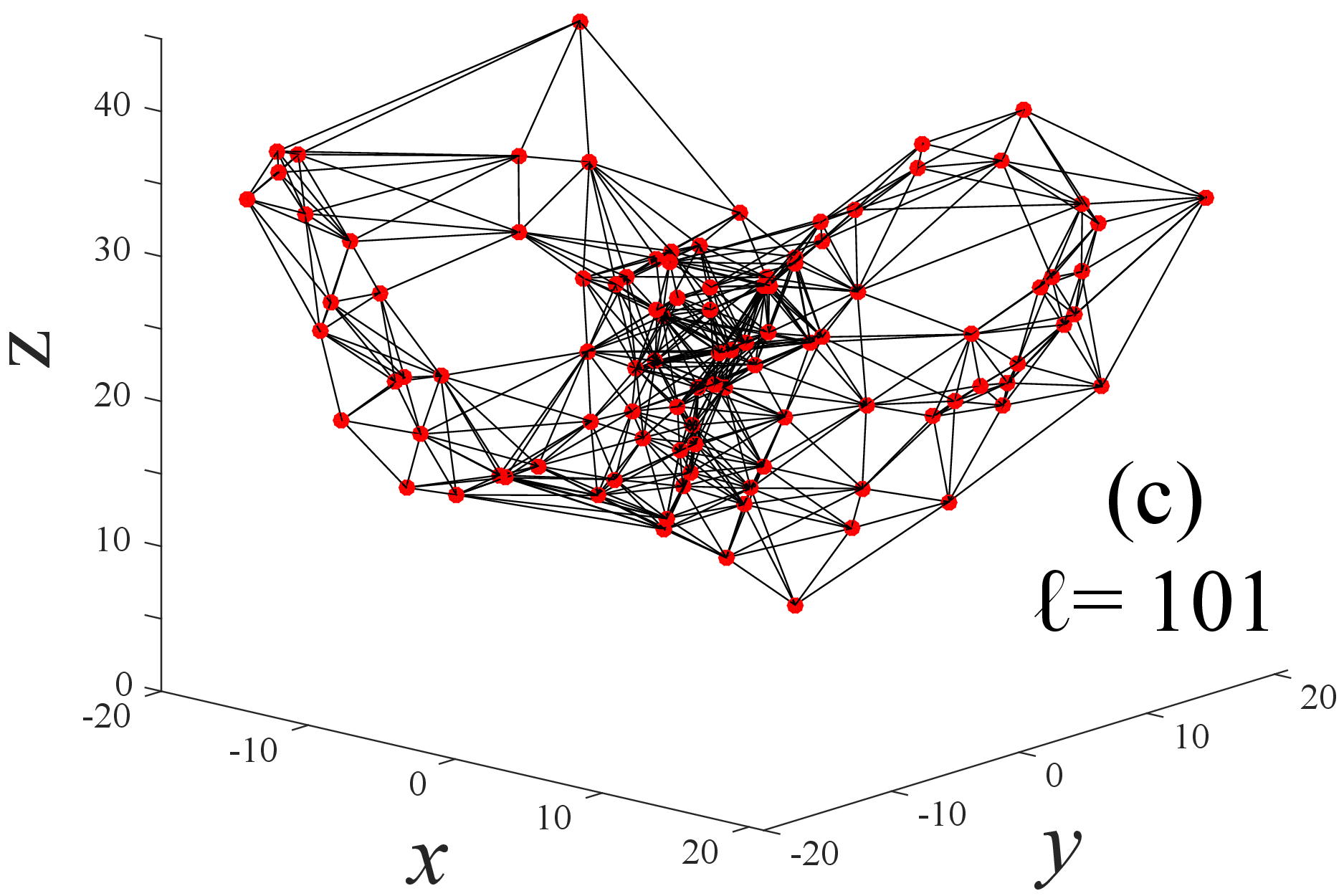}
   %             \caption{$l=101$}
        \end{subfigure}
                \begin{subfigure}[b]{0.3\textwidth}
                \includegraphics[width=\textwidth]{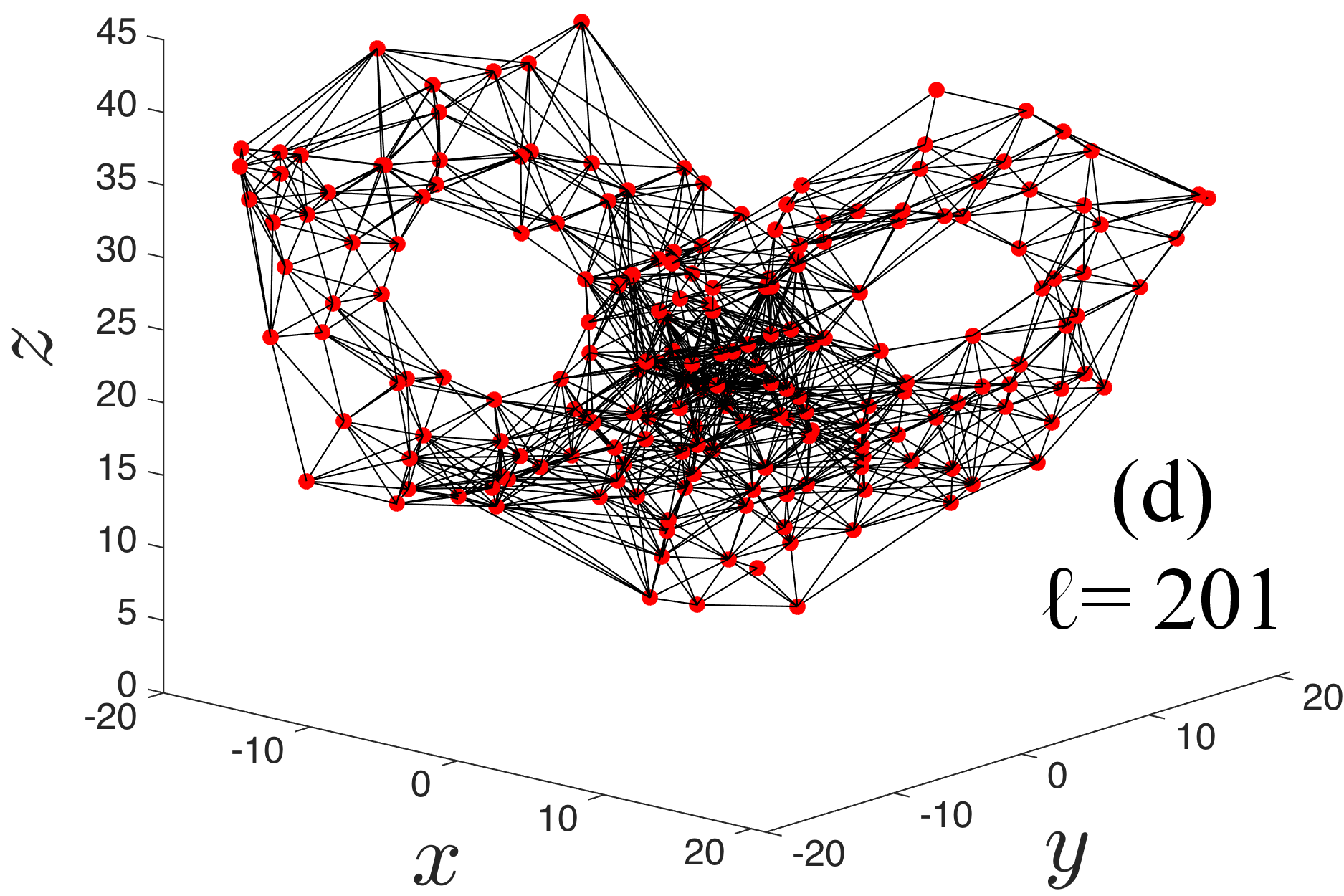}
   %            \caption{$l=201$}
        \end{subfigure}
        \begin{subfigure}[b]{0.3\textwidth}
                \includegraphics[width=\textwidth]{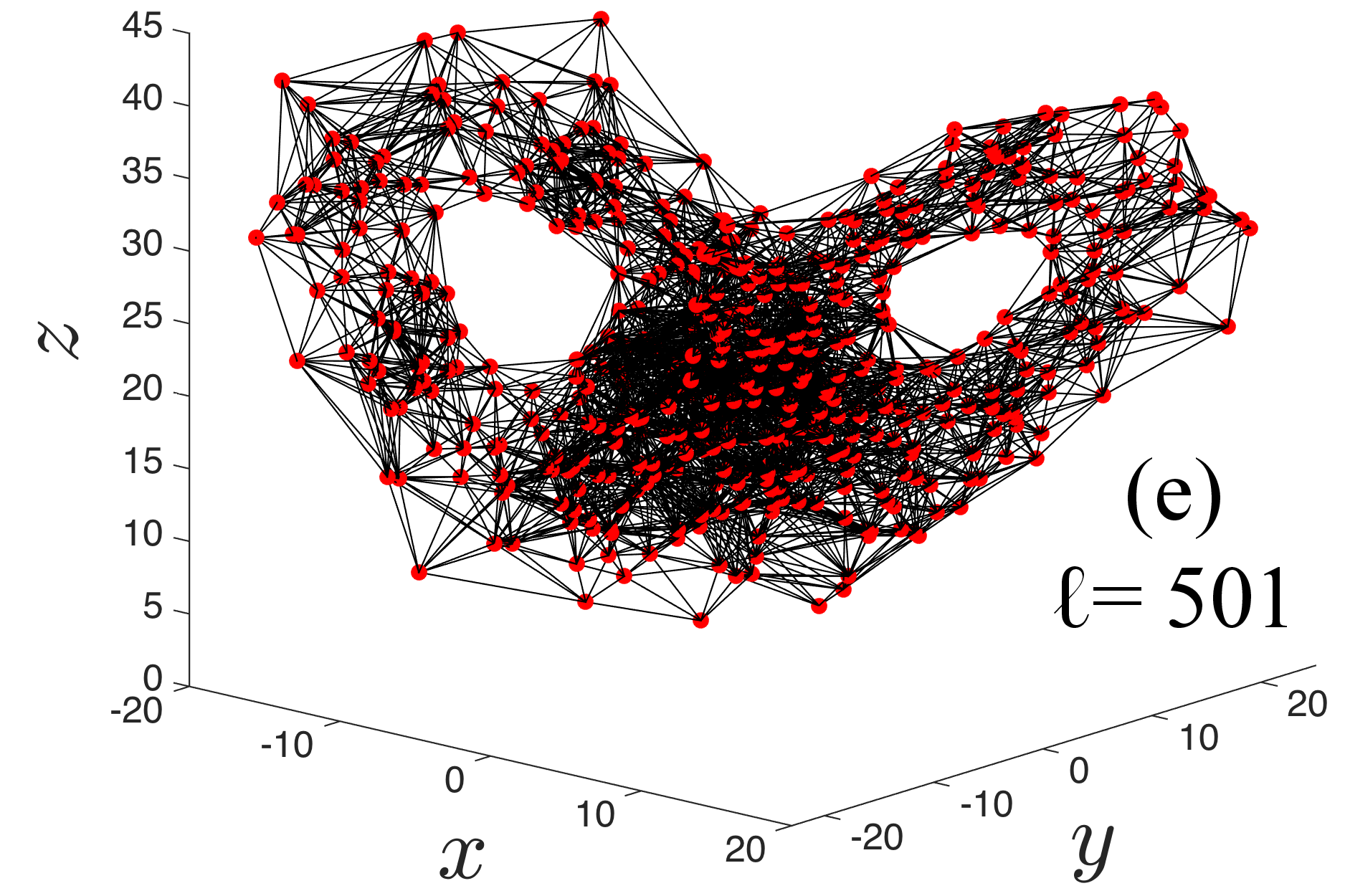}
   %             \caption{$l=501$}
        \end{subfigure}
        \begin{subfigure}[b]{0.3\textwidth}
                \includegraphics[width=\textwidth]{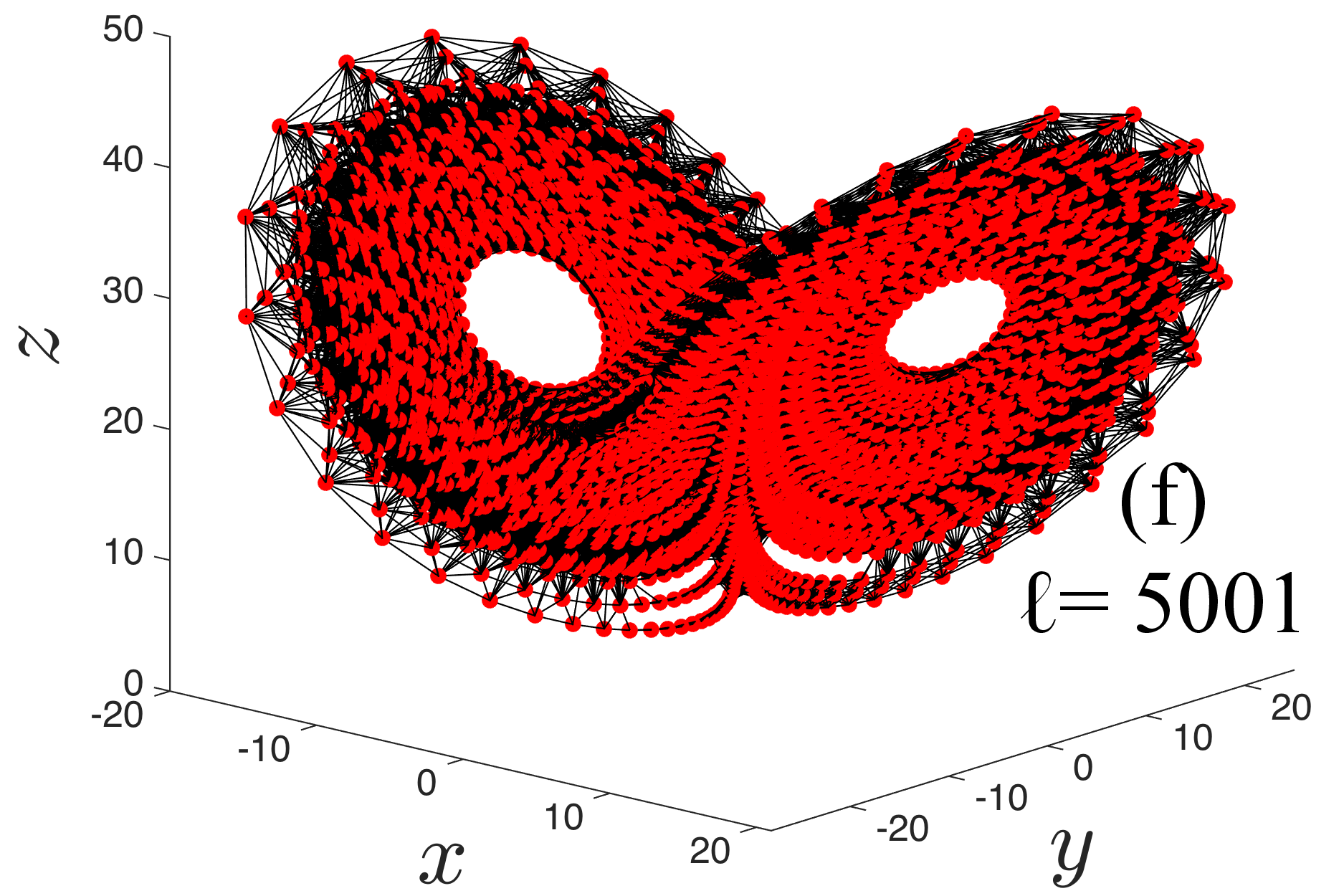}
   %             \caption{$l=5001$}
        \end{subfigure}
\caption{Varying the number of landmarks $\ell$: 
	One-skeletons of witness complexes \Eq{fuzzyComplex} constructed from
	the trajectory of \Fig{Lorenz63Example}(a) with $\eps=1.2$ and 
	six values of $\ell$.}
\label{fig:witness-complexes-different-l}
\end{figure}

In practice, we use persistence to select the number
of landmarks: i.e., given a trajectory, we build the complex for
different $\ell$ values, calculate the homology, and choose a value at
which the results stabilize.
%%%%%
In our experience, the homology of the witness complex appears to be
highly robust to the landmark spacing, but that issue will require
more exploration, as described in \Sec{concl}.
%%%
% \alert{[[make sure we actually do that -- or remove this forward pointer]]}
%%%
\Fig{witness-complexes-different-l} shows a series of examples: fuzzy
witness complexes constructed from the trajectory of
\Fig{Lorenz63Example}(a) increasing $\ell$ from $26$ to $5001$.
As one would expect, the complexes fill in as the number of landmarks
increases.  Visually, the two main holes become apparent in the
$\ell=101$ complex, which begs the question: does one really need the
extra structure of an $\ell \gg 100$ complex if the goal is to resolve
the large-scale topology of the attractor? 

%%%%%
%%%%%%

The number of landmarks required for this will obviously depend, in a
complex way, on the structure of the underlying invariant set.  In the
simplest case---if this set were a Riemannian submanifold---the
results of \cite{Nyogi08} imply that, if the ``feature size" of the
manifold is not too small, a well-defined number of sample points are
needed to reconstruct the topology of the manifold---with high
probability--- from a \v{C}ech complex with a given ball size.  It is
also possible to guarantee that the witness complex has the same
topology as an alpha complex, under appropriate conditions on the
density of the sets $L$ and $W$ and with appropriate selection of
$\alpha$ and $\eps$ \cite{Alexander15,zach-phd}.

For the case at hand, however, the invariant set is a fractal, and it
is impossible to obtain all of its structure from a finite sample,
even though we expect to resolve more of the complexity with more
landmarks.  The witness complex in
\Fig{witness-complexes-different-l}(f), for instance, resolves some of
the spiral-shaped gaps in the wings of the attractor.  For a fixed
data set, however, there is a fundamental limit: it does not make
sense to think about resolving the fine-grained structure of an object
beyond the limits that are inherent in a finite-length,
finite-precision sampling of that object.  Since all data have these
limitations, and all real-world data are noisy, it makes sense to
content ourselves with an approximate topology.  In this case, the
roughest topology corresponds to the Williams branched manifold model
of the Lorenz attractor \cite{Williams79}, which has the homology of a
figure-eight, i.e. $\beta_0 = 1$ and $\beta_1 = 2$.  This topology
corresponds to that implied by longest bars in the persistence diagram
of \Fig{Lorenz63-barcodes}(b).

The witness complex elegantly balances effectiveness and efficiency in
topological data analysis.  Its use of a small number of landmarks
sidesteps the issues that arise when one builds a fine-grained complex
that touches every data point, which is
computationally expensive and noise-sensitive.  One
could also use cubical complexes to address those issues, as in
\cite{mischaikow99}, but the witness complex is computationally less
demanding, both because it naturally follows the data and because of
the ``tight'' way that a simplicial complex covers a space.  Among
other things, the dimension of each simplex in the clique complex can be restricted
to be just high enough to cover the corresponding part of the invariant set, whereas all of the
grid elements in the cubical case necessarily have the dimension of
the ambient space.  This means that computational homology algorithms
like {\tt javaPlex} not only have fewer cells to process in the
simplicial case, but also far fewer neighbors to check---e.g., during
computations of isolating neighborhoods.  See \cite{Alexander15}[Appendix C] 
for a detailed analysis of the associated computational costs.
Witness complexes have begun
to see some use in topological data analysis \cite{Carlsson14},
but they are not completely immune to the foibles of real-world data.
Complexes constructed using the standard witness relations, for
example, can contain ``false positive'' edges due to distant witnesses
(in the case of the weak relation) and ``false negative'' edges
because of the strong relation's very stringent requirement on shared
witnesses, which is senseless in noisy, incompletely sampled data.
The fuzzy witness relation used here is intended to mitigate these
issues.  Of course, this approach is not without disadvantages.  Both
the landmark spacing $\Delta t$ that makes witness complexes
computationally efficient and the $\eps$ that makes the \emph{fuzzy}
witness complex robust with respect to noise are parameters that one
needs to choose, and choose well.  Moreover, these choices interact, as
described in the following section.

The examples presented in this section involve a full trajectory from
a dynamical system.  In the real world, however, one is generally
working with reconstructions of scalar time-series data---structures
whose topology is guaranteed to be identical to that of the underlying
dynamics if the reconstruction process is carried out properly.  But
what if the dimension $m$ does \emph{not} satisfy the requirements of
the theorems?  Can one obtain useful results about the topology of
that underlying system using the ideas, even if those dynamics are not
properly unfolded in the sense of \cite{packard80,Takens81,sauer91}?
It is to this issue that we turn next.

%%%%%%%%%%%%%%%%%%%%%%%%%%%%%%%%%%%%%%%%
\section{Topologies of reconstructions}
\label{sec:meat}

Scalar time-series data from a dynamical system are a projection of the
$d$-dimensional dynamics onto $\bR^1$---an action that does not
automatically preserve the topology of the object.  The method of
delays allows one to reconstruct the underlying dynamics, up to
diffeomorphism, if the reconstruction dimension is large enough.  
% 
% In a sense, this amounts to inverting the projection.  
% 
There are a number of conditions for the successful execution of this
procedure, as mentioned in \Sec{intro}.  According to \cite{Takens81},
seven dimensions ($m>2d$) are almost always sufficient to reconstruct
the structure of the Lorenz attractor from the time series of
\Fig{Lorenz63Example}(b); the looser bounds of \cite{sauer91},
however, suggest that $m=5$ is sufficient, since the box-counting
dimension of that attractor is 2.06.  Since the state-space dimension
is generally unknown and one needs an embedding to compute $d_{box}$,
choosing $m$ is a primary challenge in nonlinear time-series analysis.
The question we wish to address in this section is: can one use the
witness complex to obtain a useful, coarse-grained description of the
topology from lower-dimensional reconstructions---say, the basic
connectivity or number of holes in an attractor that are larger than a certain
scale?

% \cmt{Liz to do: Tie the notions of
%   ``useful'' and ``approximation;'' back to the ``what {\bf is} the
%   topology, anyway, if we have finite-length, finite-precision data?''
%   point...that it makes no sense to ask about the ``full'' topology in
%   that situation.  Use that to justify why it is both appropriate and
%   useful to try to get the large-scale features right.}

\begin{figure}[ht]
 \centering
 
         \begin{subfigure}[b]{0.45\textwidth}
                \includegraphics[width=\textwidth]{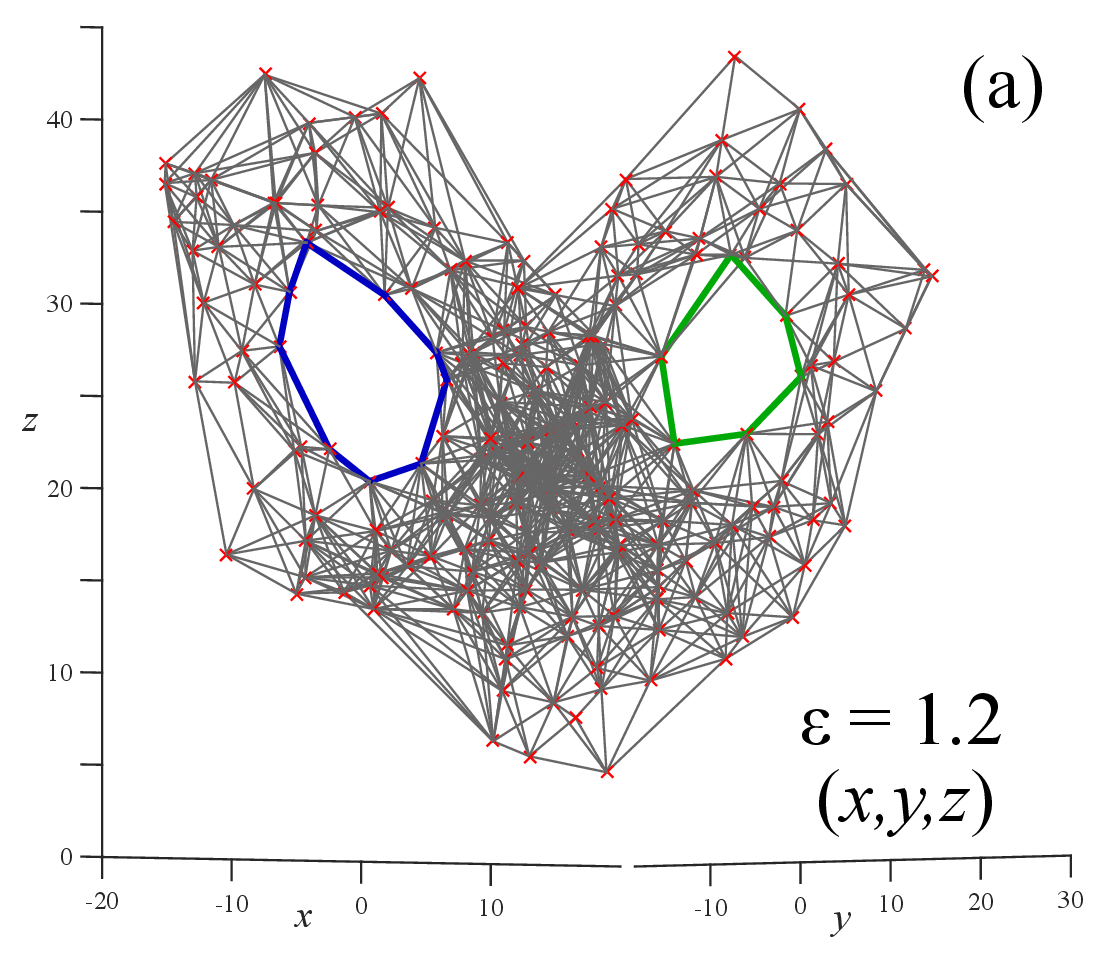}
                %\caption{Witness complex of original attractor}
        \end{subfigure}%
                 \begin{subfigure}[b]{0.45\textwidth}
                \includegraphics[width=\textwidth]{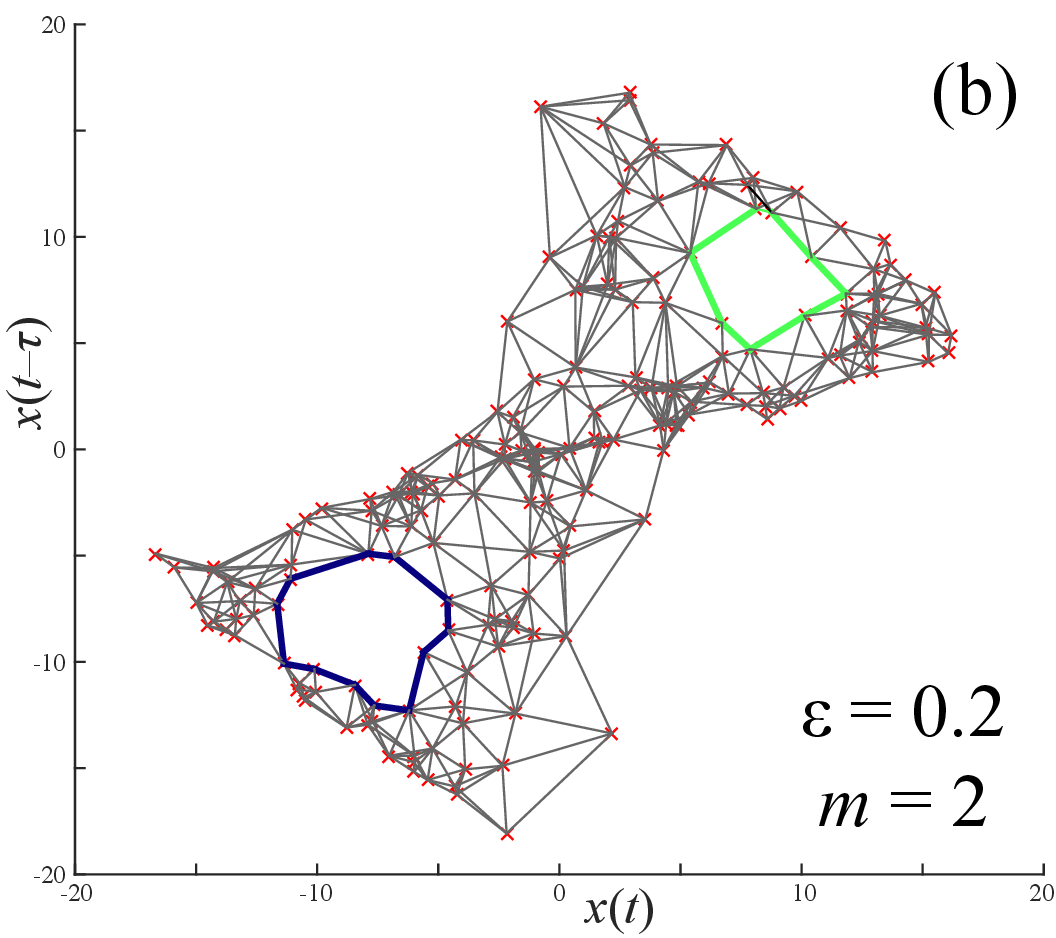}
                %\caption{Witness complex of 2-reconstruction}
        \end{subfigure}%
        \vspace*{3mm}
 
          \begin{subfigure}[b]{0.45\textwidth}
                \includegraphics[width=\textwidth]{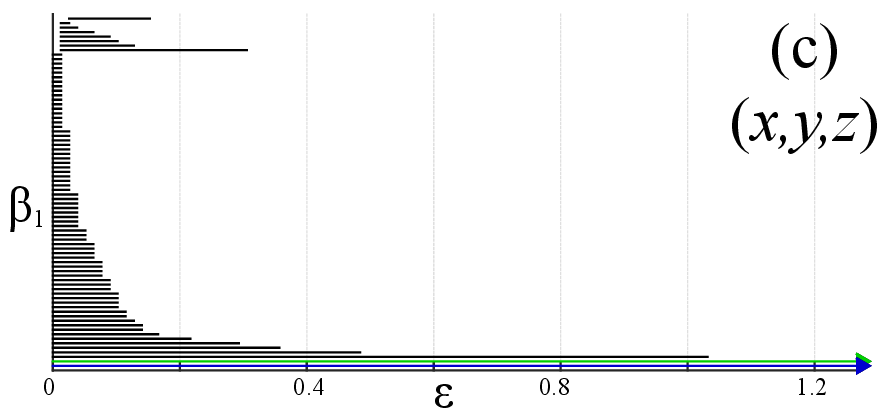}
                %\caption{$\beta_1$ barcode}
        \end{subfigure}%
                 \begin{subfigure}[b]{0.45\textwidth}
                \includegraphics[width=\textwidth]{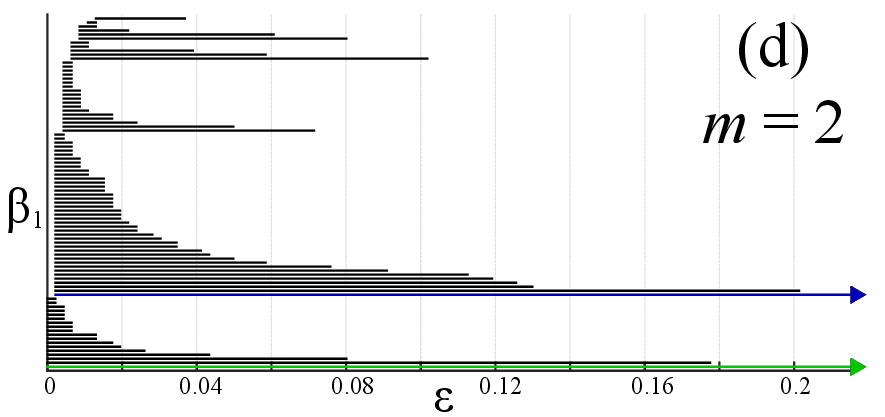}
                %\caption{$\beta_1$ barcode}
        \end{subfigure}%
        
% \mbox{\subfigure{\includegraphics[width=2.5in]{WCLorenz63Fulll200e12}
% \hspace*{5mm}
% \subfigure{\includegraphics[width=2.5in]{WCLorenz63m2l200e02}}}}
% \mbox{\subfigure{\includegraphics[width=2.5in]{Lorenz63Fulll200r432e12} 
% \hspace*{5mm}
% \subfigure{\includegraphics[width=2.5in]{Lorenz63m2l200r432e02}}}}
\caption{One-skeletons of the witness complexes (top row) and barcode
  diagrams for $\beta_1$ (bottom row) of the Lorenz system.  The plots
  in the left-hand column were computed from the three-dimensional
  $(x,y,z)$ trajectory of \Fig{Lorenz63Example}(a); those in the
  right-hand column were computed from a two-dimensional ($m=2$)
  delay-coordinate reconstruction from the $x$ coordinate of that
  trajectory with $\tau = 174 T$.  In both cases, $\ell = 201$ equally
  spaced landmarks (red {\color{red}{$\times$}}s) were used.  Both
  complexes have two persistent nonbounding cycles (green and blue
  edges) but the $2D$ reconstruction requires only $\approx1900$
  simplices to resolve those cycles (at $\eps = 0.2$), while the full
  $3D$ trajectory requires $\approx 7000$ simplices (at $\eps = 1.2$)
  to eliminate spurious loops.}
\label{fig:homology-figure1}
\end{figure}

The short answer to that question is yes.
Figure~\ref{fig:homology-figure1} shows a side-by-side comparison of
witness complexes and barcode diagrams for the Lorenz trajectory of
\Fig{Lorenz63Example}(a) and a two-dimensional reconstruction ($m=2$)
using the $x$ coordinate of that trajectory.  For the full $3D$
trajectory on the left, {\tt javaPlex} needed $6942$ simplices to
resolve the two main holes in the attractor, with $\eps=1.2$.  For the
$2D$ reconstruction on the right---constructed with $\tau = 0.174$,
the first minimum of the time-delayed mutual information---the witness
complex with $\eps = 0.2$ has only $1916$ simplices but has the same
homology as the $3D$ complex.\footnote
%%%%%
  {The difference between the
  $\eps$ values that yields a persistent result in these two cases
  makes sense because the diameters of the true and reconstructed attractors are different.
  This issue is discussed further below.}  
%%%%%
\emph{In other words, the correct large-scale homology is accessible from a witness 
complex of a $2D$ reconstruction, in a computationally efficient manner, even though 
(a) that complex does not involve all of the data points and 
(b) the reconstruction does not satisfy the conditions of the
associated theorems.}

%%%%%%%%%
%%%%%%%%

And that leads to the central question of this paper: how does the
topology of the witness complex change with the reconstruction
dimension $m$?  Intuitively,
one would think that the homology of the witness complex would change
until $m$ became large enough to correctly unfold the topology of the
underlying attractor, and then stabilize.  In practice, however,
if $m$ is too large, the so-called ``curse of
dimensionality,'' when a finite amount of data is spread over a large
volume, will destroy the fidelity of the complex. 
Additionally, the effect of noise will be amplified and
the computational expense will grow with $m$.
For all of these reasons, it would be a major advantage if one could
obtain useful information about the homology of the underlying
attractor from a low-dimensional delay-coordinate reconstruction of
scalar time-series data.

\begin{figure}
        \centering
        \begin{subfigure}[b]{0.45\textwidth}
                \includegraphics[width=\textwidth]{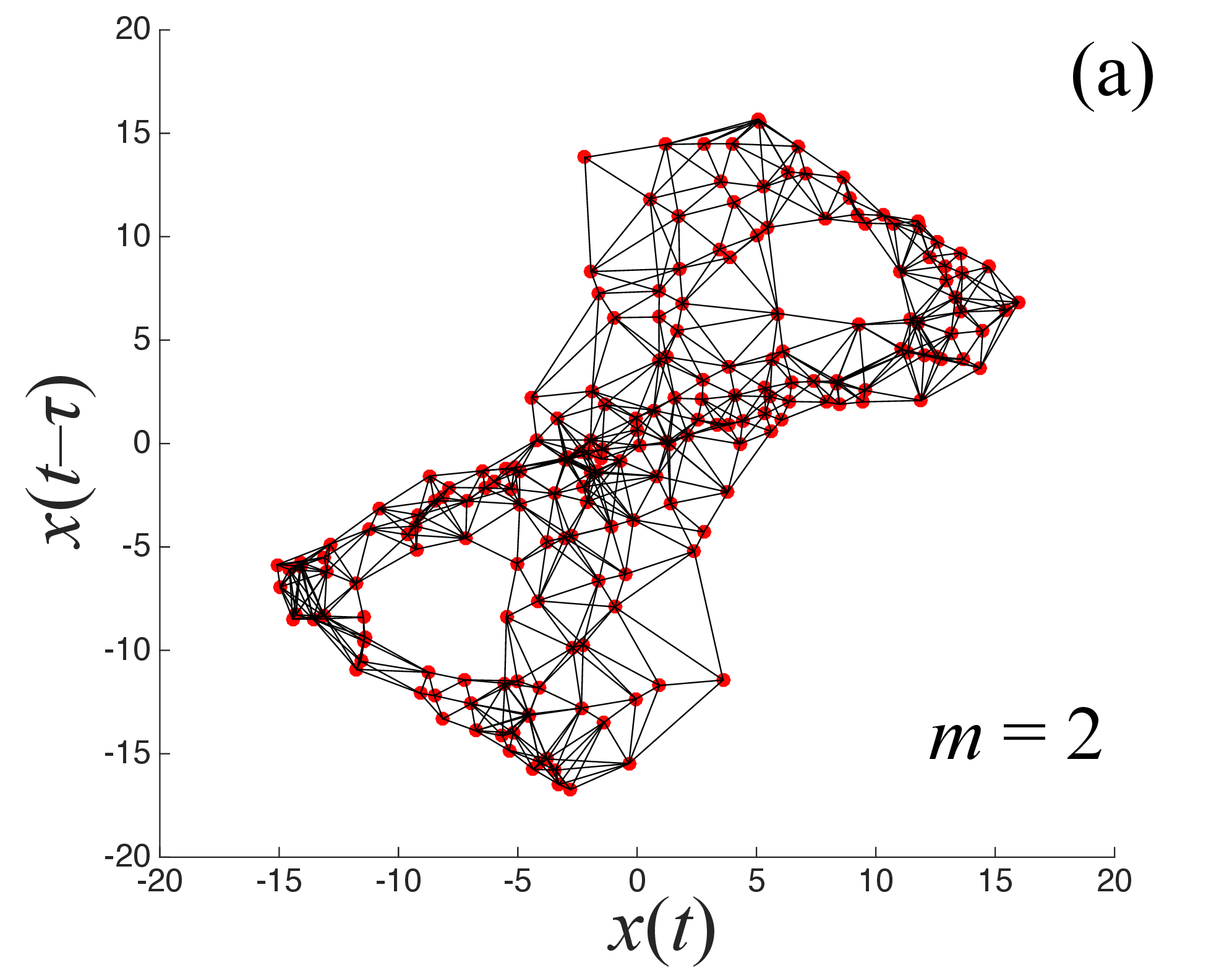}
        \end{subfigure}%
        \begin{subfigure}[b]{0.45\textwidth}
                \includegraphics[width=\textwidth]{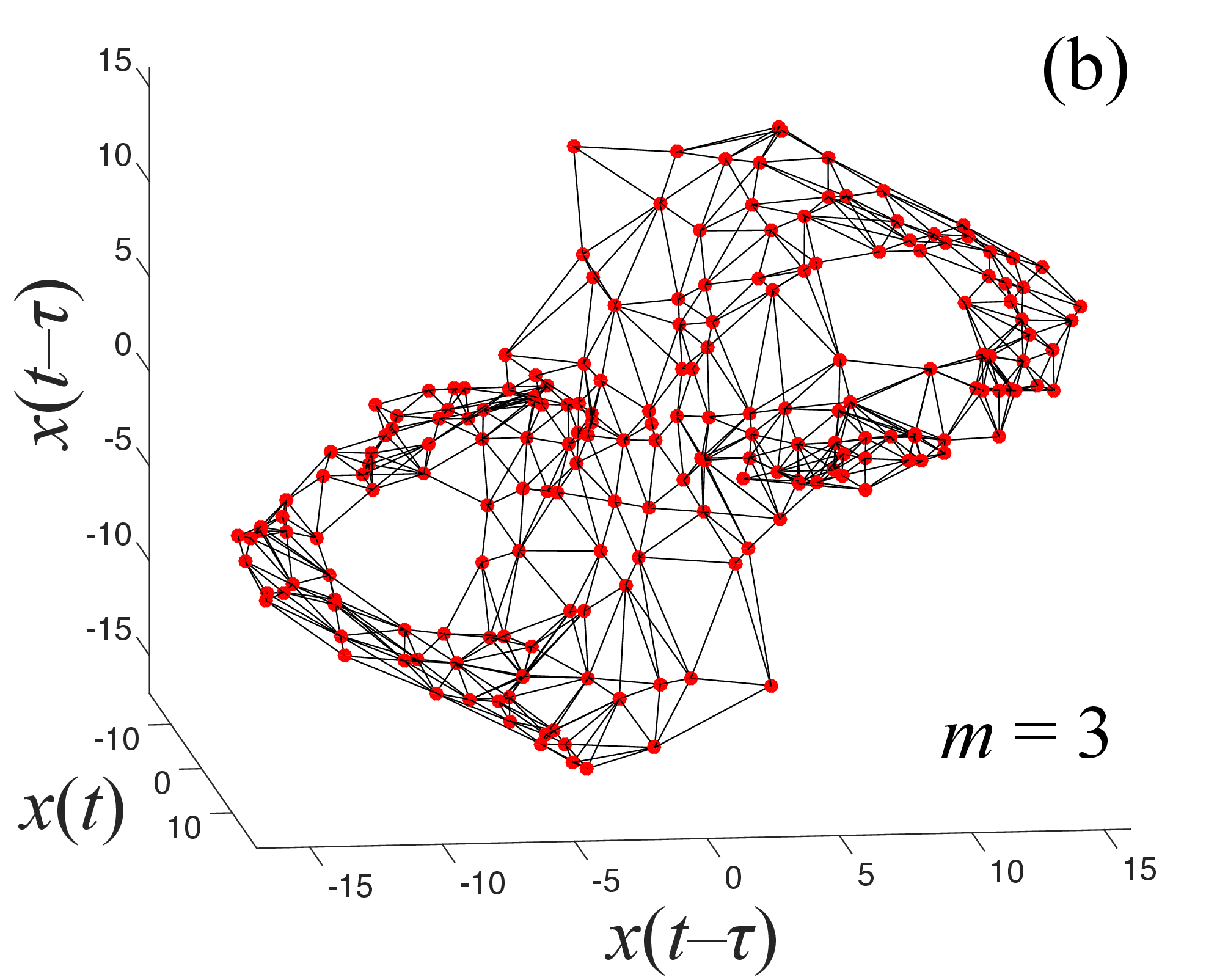}
        \end{subfigure}
\caption{The effect of reconstruction dimension: One-skeletons of
  witness complexes of different reconstructions of the scalar time
  series of \Fig{Lorenz63Example}(b).  Both reconstructions use $\tau
  =0.174$, the first minimum of the average time-delayed mutual
  information\cite{fraser-swinney}, $\ell=198$ equally spaced
  landmarks (red dots), and $\xi = 0.54\%$, as defined in \Eq{scaledXiDef}.}
  %$\eps$ values are scaled as explained in the text.
\label{fig:LorenzEmbedSkeletons}
\end{figure}

Again, it appears that this is possible.
Figure~\ref{fig:LorenzEmbedSkeletons} shows witness complexes for
$m=2$ and $m=3$ reconstructions of the Lorenz time series of
\Fig{Lorenz63Example}(b).
The barcodes for first two Betti numbers for these two complexes, as computed using
{\tt javaPlex}, have similar structure: the complexes become connected ($\beta_0 = 1$) at a small value of $\eps$, and the dominant, persistent homology corresponds to the two primary holes ($\beta_1 = 2$) in the attractor. 
Note, by the way, that \Fig{LorenzEmbedSkeletons}(a) is not simply a
$2D$ projection of \Fig{LorenzEmbedSkeletons}(b); the edges in each
complex reflect the geometry of the witness relationships in different
spaces, and so may differ.
(The \emph{landmark} set in the $2D$
reconstruction \emph{is} a projection of the landmark set in the $3D$
reconstruction, however, because we apply the same ``lifting''
operation to all points.)
%%%%
Higher-dimensional reconstructions---not
easily displayed---have the same homology for suitable choices of
$\eps$, though for $m >5$, it is necessary to increase the number of
landmarks to obtain a persistent $\beta_1 = 2$.

%%%%%%%%%%%
%%%%%%%%%%

That brings up an important point: if one wants to sensibly compare
witness complexes constructed from different reconstructions of a
single data set, one has to think carefully about the $\ell$ and
$\eps$ parameters.  Here, we used persistence to choose a good value
of $\ell$.  We found that the results were robust with respect to
changes in that value, across all reconstruction dimension values in
this study, so we fix $\ell \approx 200$ for all the experiments
reported here.\footnote
%%%
{The precise value varies slightly because the length of a
  trajectory reconstructed from a fixed-length data set decreases with
  increasing $m$ (since one needs a full span of $m \times (\tau/T)$ data
  points to construct a point in the reconstruction space).}
%%%
Since the number of data points required to properly sample an object should
generally grow with dimension \cite{tsonisdatabound}, this will
require more exploration, as mentioned in \Sec{concl}.

In the experiments in the previous section, the scale parameter $\eps$
was given in absolute units.  To generalize this approach across
different examples, it makes sense to compare reconstructions with
$\eps$ chosen to be a fixed fraction,
\beq{scaledXiDef} \eps =  \xi \, \diam{W} \eeq
of the diameter, $\diam{W}$, of the set $W$.  For example, for the
full $3D$ attractor in \Fig{Lorenz63Example}(a),
\[
	\diam{W_{xyz}} = \sqrt{(x_{max} - x_{min})^2 +
  (y_{max} - y_{min})^2 + (z_{max} - z_{min})^2} = 75.3 ,
\]
%%%%%% lorenz63 distances/diamaters
%dx = max(x) - min(x) = 37.0235
%dy = max(y)-min(y) = 50.2910
%dz = max(z)-min(z) = 42.0139
%dA = sqrt( (dx)^2 +(dy)^2+(dz )^2) = 75.2668
%%%%%
so the $\eps$ values used in \Fig{Lorenz63-barcodes}---$0.017 \le \eps
\le 1.7$ in absolute units---translate to $2.3 \times 10^{-4} \le
\xi \le 0.023$ in this diameter-scaled measure.

The diameter of the reconstruction varies in a natural way with the dimension $m$.
Since delay-coordinate reconstruction of scalar data unfolds
the full range of those data along every added dimension, the diameter
of an $m$-dimensional reconstruction will be
\[
	\diam{W_m} = \sqrt{m(x_{max}-x_{min})^2} = 37.0 \sqrt{m} ,
\]
where $x$ represents the scalar time-series data.
%If one had perfect data,
%that would not call for any change in the scale parameter of the
%analysis.  We have only a finite sample, however---in length,
%precision, and spacing---so we 
Since this unfolding will change the geometry of the reconstruction,
we need to scale $\eps$ accordingly.  The witness
complexes in \Fig{LorenzEmbedSkeletons} were constructed with a fixed
value of $\xi = 0.54\%$.
%
% in \Eq{scaledXiDef}. 
%
Thus, for \Fig{LorenzEmbedSkeletons}(a), $\eps = 37.0\sqrt{2}(0.0054)
= 0.283$ in absolute units, while for \Fig{Lorenz63Example}(b),
$\diam{W_3}=37.0\sqrt{3}$ and $\eps= 0.346$.
Together with the even temporal spacing of landmarks, this scaling of
$\eps$---which is used throughout the rest of this section---should
allow the witness complex to adapt appropriately to the effects of
changing reconstruction dimension and finite data.
%%%%%
% \alert{Need to talk about that and figure out how say it better and
%   more formally, obviously, and probably wave hands about future work
%   regarding whatever we say.  And we do need to mention the
%   finite-data issue in there somewhere, since we wouldn't have to
%   scale \eps if the data were perfect.}
%%%%

\begin{figure}[ht!]
        \centering
        \centering
        \begin{subfigure}[b]{0.45\textwidth}
		\includegraphics[height=1.3truein]{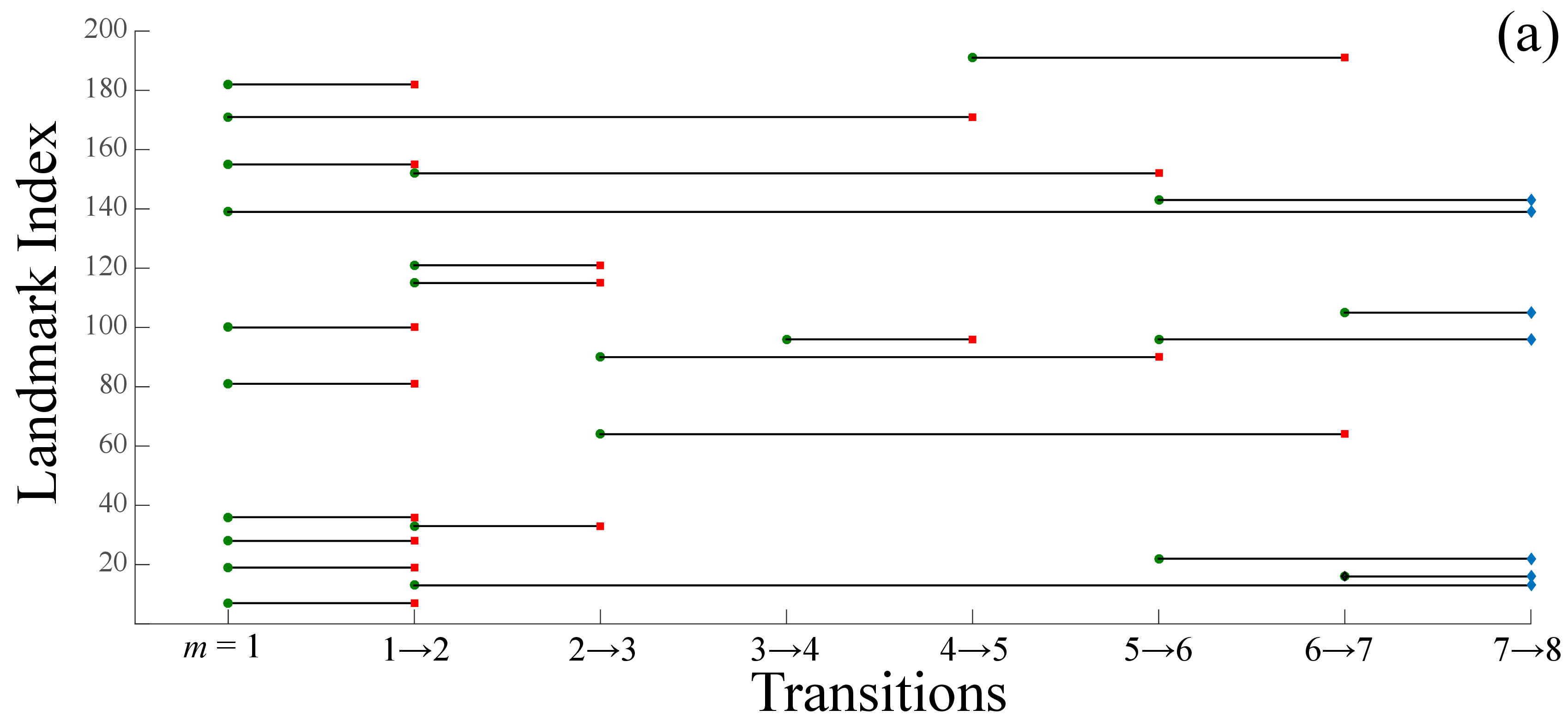}

		%\caption{}
        \end{subfigure}%
\hspace*{25mm}
        \begin{subfigure}[b]{0.35\textwidth}
		\includegraphics[height=1.3truein]{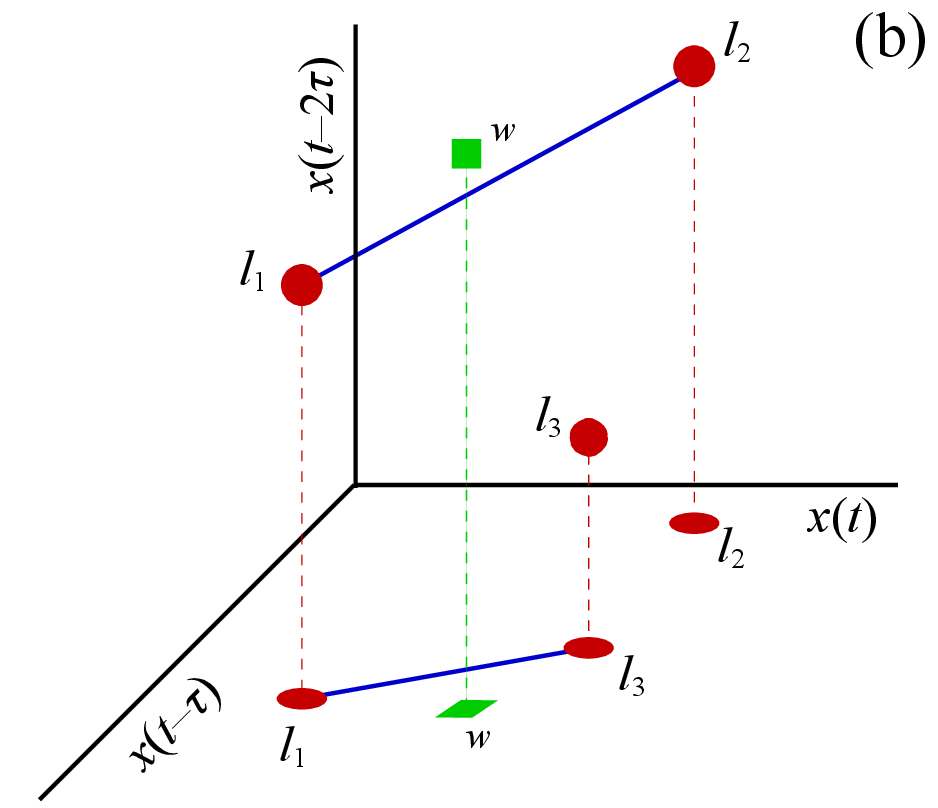}

		%\caption{}
        \end{subfigure}
        \caption{(a) Dimension barcode for edges in the witness complex of the
  reconstructed scalar time series of \Fig{Lorenz63Example}(b) that
  involve $l_0$, the first landmark, for reconstructions with
  $m=1,\ldots, 8$.  The vertical axis is labeled with the indices of
  the remaining $197$ landmarks in the complex; a (green) circle at the
  $m-1\rightarrow m$ tickmark on the horizontal axis indicates the
  transition at which an edge between $l_0$ and $l_i$ is born; a (red)
  square indicates the transition at which that edge vanishes from the
  complex.  A (blue) diamond at the right-hand edge of the plot
  indicates an edge that was still stable when the algorithm
  completed. For all reconstructions, $\tau =0.174$, $\ell=198$, and
  $\xi=0.54\%$.  (b) Sketch of the birth of an edge at the
  $m=2\rightarrow3$ transition.}
\label{fig:EdgeStability}
\end{figure}

To formalize the exploration of the reconstruction homology and extend
that study across multiple dimensions, one can use a variant of the
classic barcode diagram that shows, for each simplex, the
reconstruction dimension values at which it appears in and vanishes
from the complex.  \Fig{EdgeStability}(a) shows such a plot for edges
that involve $l_0$, the first landmark on the reconstructed
trajectory.  A number of interesting features are apparent in this
image.  Unsurprisingly, most of the one-simplices that exist in the
$m=1$ witness complex---many of which are likely due to the strong
effects of the projection of the underlying $\bR^d$ trajectory onto
$\bR^1$---vanish when one moves to $m=2$.  There are other short-lived
edges in the complex as well: e.g., the edge from $l_0$ to $l_{120}$
that is born at $m=2$ and dies at $m=3$.  The sketch in
\Fig{EdgeStability}(b) demonstrates how edges can be born as $m$
increases: in $m=2$, $\ell_1$ and $\ell_3$ share a witness (the green
square); when one moves to $m=3$, spreading all of the points out
along the added dimension, that witness is moved far from
$\ell_3$---and into the shared region between $\ell_1$ and $\ell_2$.
There are also long-lived edges in the complex of
\Fig{EdgeStability}(a).  The one between $l_0$ and $l_{140}$ that
persists from $m=1$ to $m=8$ is particularly interesting: this pair of
landmarks has shared witnesses in the scalar data \emph{and in all
  reconstructions}.  Possible causes for this are explored in more
depth below.  All of these effects depend on $\xi$, of course;
decreasing $\xi$ will decrease both the number and average length of
the edge persistence bars.

%%%%
%%%%%%

\begin{figure}[ht!]
        \centering
        \begin{subfigure}[b]{0.6\textwidth}
                \includegraphics[width=\textwidth]{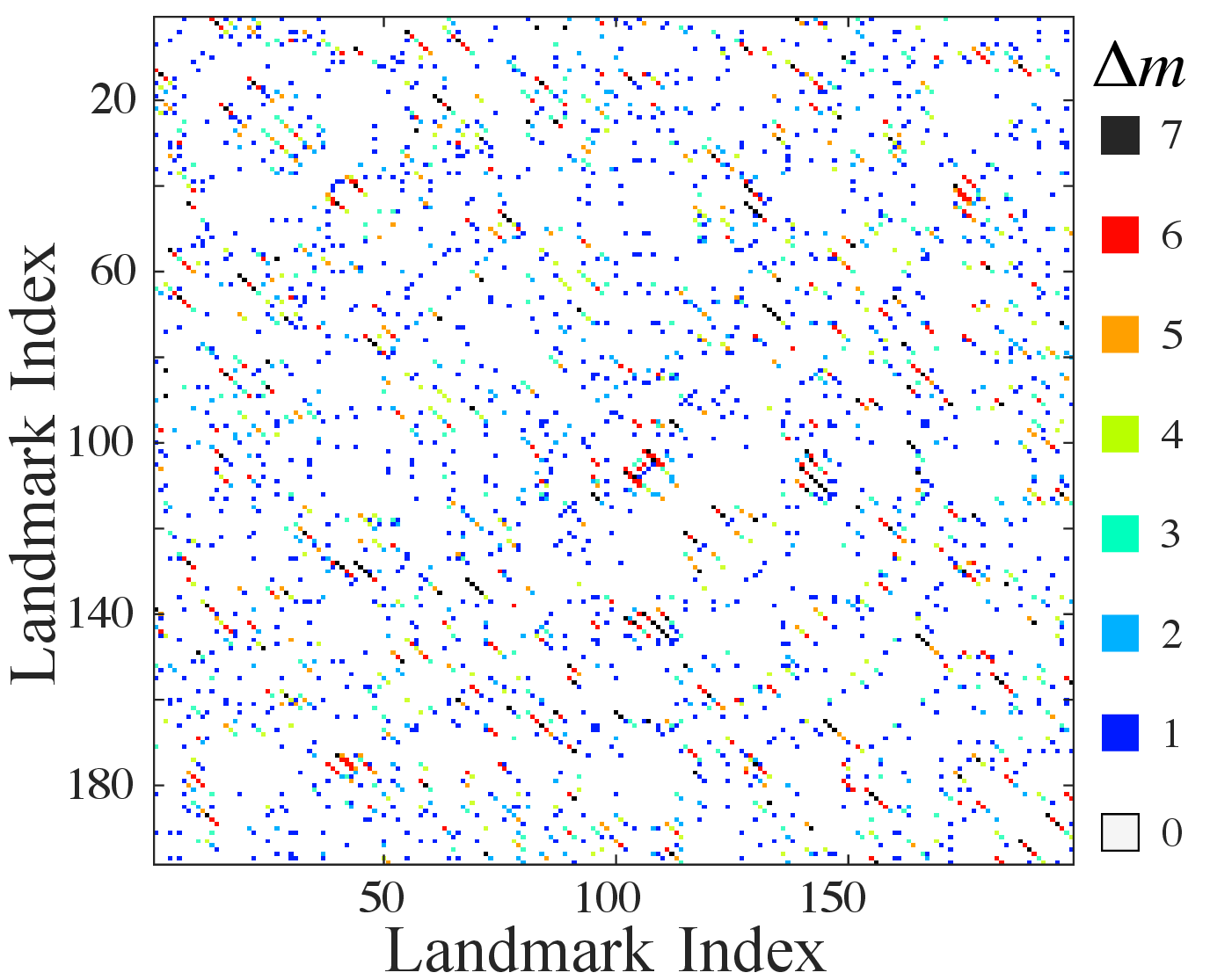}
        \end{subfigure}
% 
%         \begin{subfigure}[b]{0.45\textwidth}
%                 \includegraphics[width=\textwidth]{persistHistoScaledEps}
%         \end{subfigure}
% \includegraphics[width=0.5\textwidth]{persistDistr}
% \includegraphics[width=0.7\textwidth]{heatmapedgepersistdiagram}
\caption{Edge lifespan diagram: pixel $i,j$ on this image is
  color-coded according to the maximum range $\Delta m$ of dimension
  for which an edge exists between landmarks $l_i$ and $l_j$ in the
  witness complex of the reconstructed scalar time series of
  \Fig{Lorenz63Example}(b) for $m=1,\ldots, 8$.
% 
% Image (b) shows a histogram of how many edges in this series of
% complexes are $\Delta m$-persistent for different values of $\Delta
% m$.
% 
For all reconstructions, $\tau =174$, $\ell=198$, and $\xi=0.54\%$.}
\label{fig:BirthDeathLorenz}
\end{figure}

While this $\Delta m$ barcode image is interesting, the amount of
detail that it contains makes it somewhat unwieldy.  To study the
$m$-persistence of all of $\ell \times \ell$ edges in a witness
complex, one would need to examine $\ell$ of these plots---or condense them into a
single plot with $\ell^2$ entries on the vertical axis.  Instead, one can plot
what we call an \emph{edge lifespan diagram}: an $\ell \times \ell$
matrix whose $(i,j)^{th}$ pixel is colored according to the maximum
range of $m$ for which an edge exists in the complex between the $i^{th}$
and $j^{th}$ landmarks; see  \Fig{BirthDeathLorenz}.
If the edge $\{l_i,l_j\}$ existed in the complex for
$2 \le m<3$ and $5 \le m<8$, for instance, $\Delta m$ would be three and the
$i,j^{th}$ pixel would be coded in cyan. Edges that do not exist for
any dimension are coded white.  

%%%%%%
%%%%%%

A prominent feature of \Fig{BirthDeathLorenz} is a large number
($683$) of edges with a lifespan $1$ (blue). Of these edges, $463$
exist for $m=1$, but not for $m=2$, and thus reflect the anomalous
behavior of projecting a $2.06$ dimensional object onto a line.  This
was also seen, as described above, in the barcode of \Fig{EdgeStability}.

Another feature observed in the lifespan diagram is a number of
diagonal line segments. Note that the \emph{color} of the pixels in
the segments varies, though most of these segments correspond to edges
with longer lifespans.  These segments indicate the existence of
$\Delta m$-persistent edges $\{l_i,l_j\}, \{l_{i+1},l_{j+1}\},
\{l_{i+2}, l_{j+2}\} \ldots$.  This is likely due to the continuity of
the dynamics \cite{Alexander12}.  Recall that the landmarks are evenly
spaced in time, so $l_{i+1}$ is the $\Delta t$-forward image of $l_i$.
Thus a diagonal segment may indicate that the $\Delta t$-forward
images of (at least one) witness that is shared between $l_i$ and
$l_j$ is shared between $l_{i+1}$ and $l_{j+1}$, and so on.  The
lengths of the longer line segments suggest that that continuity fails
after $5$-$10$ $\Delta t$ steps, probably because of the positive
Lyapunov exponents on the attractor. A simple check on this reasoning,
would be an edge lifespan diagram for a dynamical system with a limit
cycle.
% Lorenz with ($a=16, \, r=500, \, b=4$).
The structure of this plot (not shown) is dominated by diagonal lines
of high $\Delta m$-persistence, with a few other scattered
one-persistent edges.  One can capture the underlying dynamical
information that gives rise to these effects more fully using what we
call the witness \emph{map}, as mentioned briefly in \Sec{concl} and
described at more length in \cite{Alexander15}.
%%%
% \alert{[[make sure we actually do that -- or remove this forward pointer]]}
%%%
% Part~(b) of the Figure shows a histogram of the number of edges that
% are $\Delta m$-persistent, as a function of $\Delta m$.  All 14 of the
% one-persistent edges in \Fig{EdgeStability}, for instance, contribute to
% the leftmost bar in this plot.  \cmt{Say something sagacious about the
%   shape of the histogram.}

\begin{figure}
        \centering
                \includegraphics[width=0.5\textwidth]{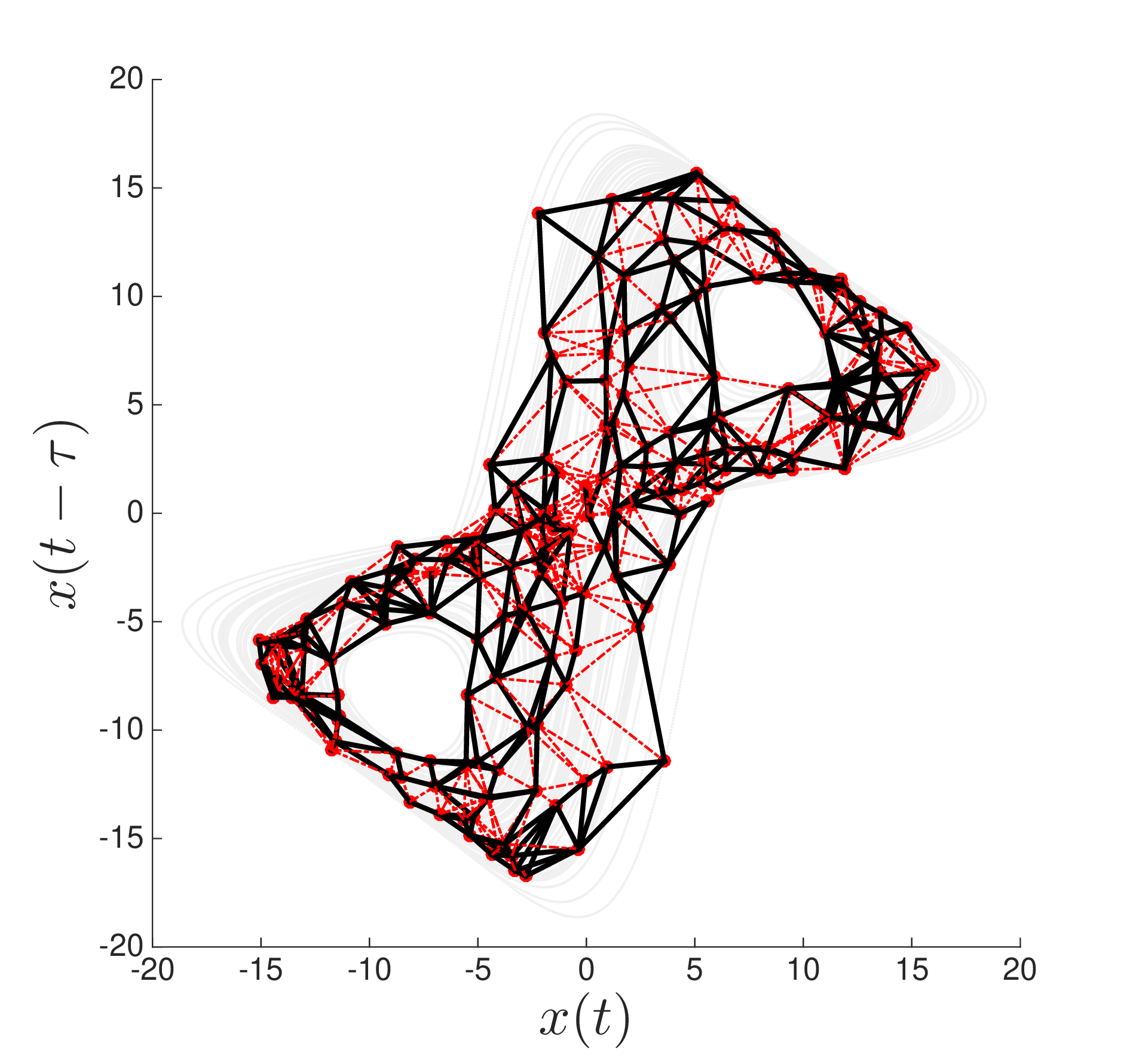}
 \caption{Witness complex of \Fig{LorenzEmbedSkeletons}(a) with $\Delta m \ge 2$-persistent edges shown as thick (black) lines, and the $\Delta m = 1$ edges as (red) dashed lines.}
%   $\Delta m=3$}
 \label{fig:delta-m}
 \end{figure}

The rationale behind studying the \emph{maximal} $m$-lifespan goes
back to one of the basic premises of persistence: that features that
persist for a wide range of parameter values are in some sense
meaningful.  To explore this, \Fig{delta-m} shows the witness complex
of \Fig{LorenzEmbedSkeletons}(a) with the $\Delta m \ge 2$-persistent
edges drawn as thicker lines: that is, edges that exist at $m=2$ and
persist at least to $m=4$.  There exists a fundamental core to the
complex that persists as the dimension grows and thus is robust to
geometric distortion, but there are also short-lived edges that fill
in the complex in accord with the local geometric structure of the
reconstruction.  Indeed, when $m=2$, the projection artificially
compresses near the origin; small simplicies fill in this region due
to the landmark clustering there. However, in the transition to
$m=3$---viz., \Fig{LorenzEmbedSkeletons}(b)---this region stretches
away from the origin, spreading the landmarks out.  There is a similar
cluster of fragile edges near the lower left corner of the complex.

This geometric evolution with increasing reconstruction dimension
leads to the death of many local edges.  Even so, the large-scale
topology is correct in both complexes of \Fig{LorenzEmbedSkeletons},
although the fine-scale topology is resolved differently by the
dimension-dependent geometry.  So while the edges with longer lifespan
are indeed more important to the core structure, the short-lived edges
are also important because they allow the complex to adapt to the
geometric evolution of the attractor and fill in the details of the
skeleton that are necessary and meaningful in that dimension.

In the spirit of the false near-neighbor method \cite{KBA92}, one
might be tempted to take the short-lived edges as an indication that
the reconstruction dimension is inadequate.  However, one computes
homology \emph{from the overall complex}.  As the example above shows,
homology is relatively robust with respect to individual edges.  The
moral of this story is that the lifespan of an edge is not necessarily
an obvious indication of its importance to the homology of the
complex; $\Delta m$-persistence plays a different role here than the
abscissa of traditional barcode persistence plots.
%The persistent edges form the foundation of the topology and the short-term edges are the 
%building blocks that allow the topology to be filled in adaptively in accordance with how the 
%geometry of the dynamics is evolving. 

%%%%%%%%%
%%%%%%%%%
% \begin{figure}
%         \centering
%              \begin{subfigure}[b]{0.49\textwidth}
%                 \includegraphics[width=\textwidth]{lorenz63ComplexEE1LS}
%                  \caption{Born at 1 died at 2 \cmt{(for $m=1$ project to y-axis!!!!)}}
%         \end{subfigure}
%         \begin{subfigure}[b]{0.49\textwidth}
%                 \includegraphics[width=\textwidth]{lorenz63ComplexG5LS}
%                  \caption{$LS\ge5$}
%         \end{subfigure}
% 
%          \begin{subfigure}[b]{0.49\textwidth}
%                 \includegraphics[width=\textwidth]{lorenz63ComplexG6LS}
%                  \caption{$LS\ge6$}
%                         \end{subfigure}
%         \begin{subfigure}[b]{0.5\textwidth}
%                 \includegraphics[width=\textwidth]{lorenz63ComplexG7LS}
%                  \caption{$LS\ge7$}
%         \end{subfigure}
%                \caption{}
% \caption{Witness complexes of $\Delta m$-persistent edges for
%   different values of $\Delta m$}
% \label{fig:delta-m}
% \end{figure}

\begin{figure}
        \centering
        \begin{subfigure}[b]{0.45\textwidth}
                \includegraphics[width=\textwidth]{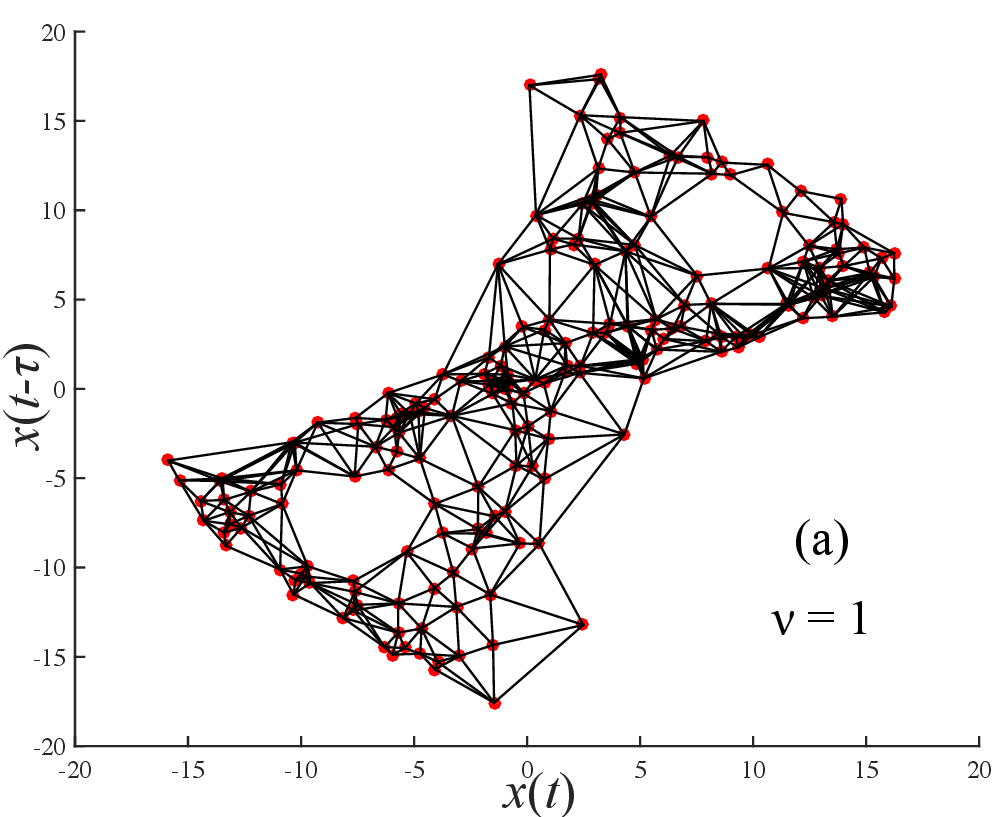}
             %  \label{fig:2dskel}
        \end{subfigure}%
        \begin{subfigure}[b]{0.45\textwidth}
                \includegraphics[width=\textwidth]{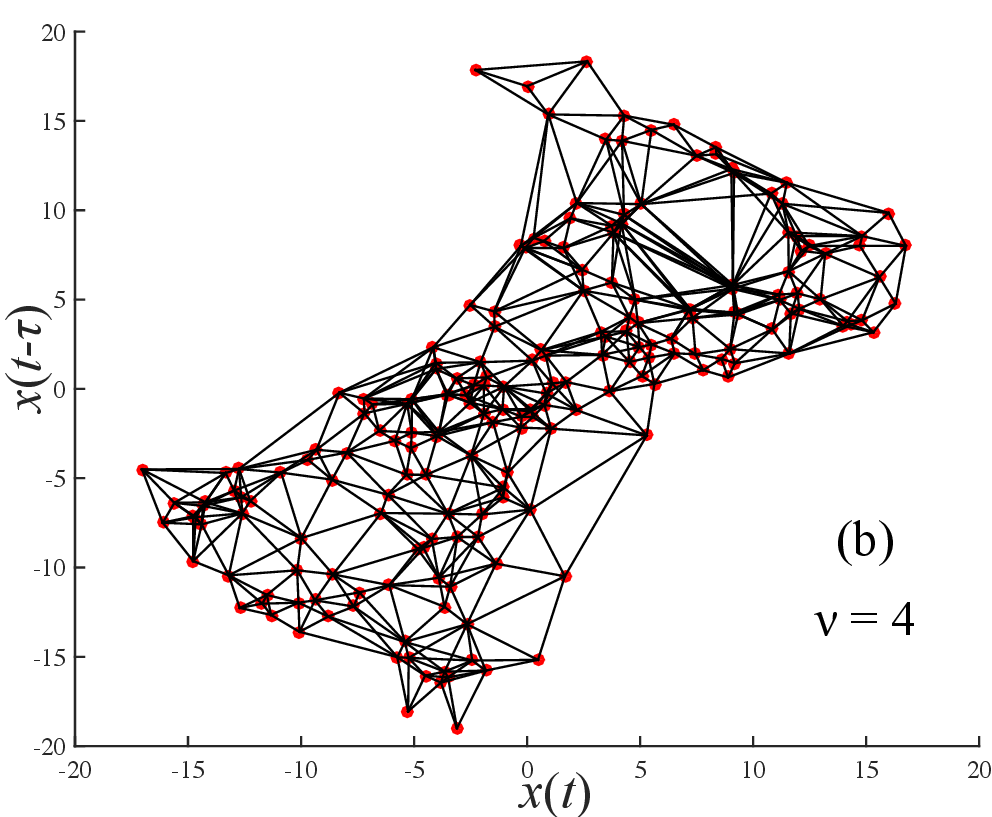}
              %  \label{fig:3dskel}
        \end{subfigure}

%\vspace*{3mm}

        \begin{subfigure}[b]{0.45\textwidth}
                \includegraphics[width=\textwidth]{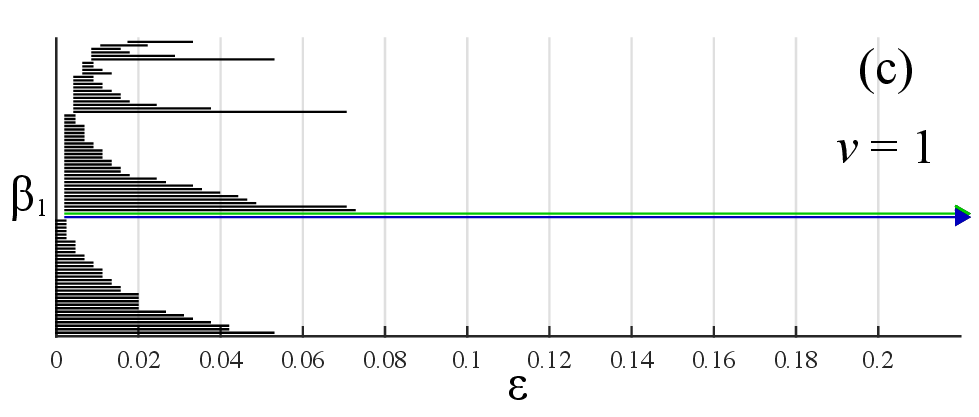}
        \end{subfigure}
        \begin{subfigure}[b]{0.45\textwidth}
			\includegraphics[width=\textwidth]{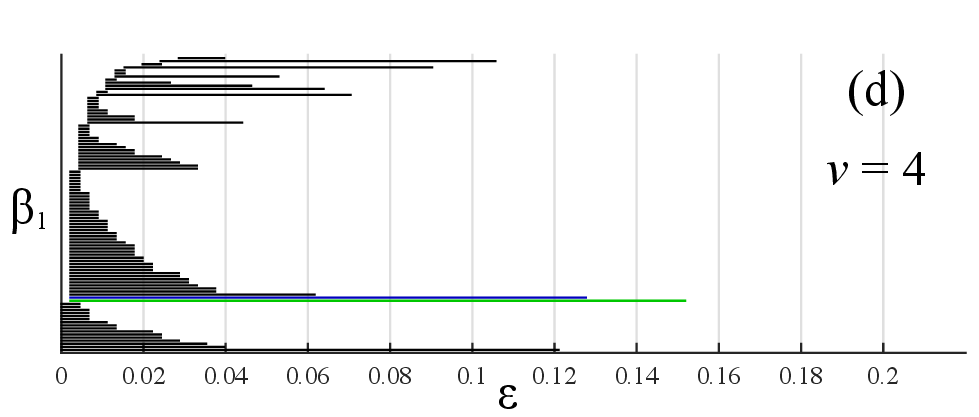}
        \end{subfigure}		
\caption{The effect of noise: One-skeletons and barcodes of witness
  complexes for $m=2$, $\xi=0.41\%$, $\tau=174 T$, $\ell=201$
  reconstructions of the scalar time series of
  \Fig{Lorenz63Example}(b) with added uniform noise of size $\nu = 1$
  (a and c) and $\nu = 4$ (b and d)---respectively, 1.9\% and 7.6\% of the
diameter of the reconstructed attractor.}
\label{fig:noise}
\end{figure}

A closely related issue is noise, which is always present in
real-world data and can disturb the geometric relationships between
points in the complex.  To study the effect on the fuzzy witness
complex, we add uniformly distributed noise on the interval
$[-\nu/2,\nu/2]$ to each point of the trajectory of
\Fig{Lorenz63Example}(b), and then perform a delay-coordinate
reconstruction using $m=2$.  Visually, \Fig{noise}(a), where $\nu =1$,
shows a surprisingly similar complex to the noise-free case of
\Fig{homology-figure1}(b); however, the structure is clearly different
in \Fig{noise}(b), where $\nu=4$.  Comparing the barcodes, one sees
that the fine-scale structure is, unsurprisingly, washed out for
smaller values of $\eps$ than in the noise-free case.  In particular,
the persistence intervals of small-scale loops are decreased.  It is
encouraging, however, to see that the two major holes persist over a
wide range of $\eps$ even when the noise level is close to 2\% of the
diameter of the attractor.  The larger noise level (7.6\% of that
diameter) is enough to destroy the two large loops so that the
complex becomes acyclic when $\eps > 0.152$.

%%%%%%%
%%%%%%

The relative immunity to noise is a general feature of persistent
homology for point cloud data \cite{Ghrist08}.  The additional
robustness of the witness complex has two sources: the fuzziness
parameter and the fact that multiple points can witness a given
simplex.  Even if noise moves one witness out of the shared region
that is defined by $\eps$, it may not move \emph{all} of the witnesses
to a particular simplex out of that shared region.  This suggests a
possible noise mitigation technique: build a complex that only
contains simplices that have at least $n$ witnesses---and perhaps
explore $n$-persistence of the associated homology.  We plan to
explore these ideas in future work.

Assessing the change in topology with changing reconstruction
dimension is a new flavor of persistence---an idea that has
traditionally been applied in the context of scale parameters like
our fuzziness parameter $\eps$.  Recall
that persistent homology is based on the idea of a filtration, i.e., a
nested set of complexes.  Most of the standard simplicial complex
constructions for data have a natural parameter that gives rise to a filtration;
for example, the radius for a sequence of \v{C}ech complexes.  This
filtration can be used to define persistent homology
\cite{Zomorodian05} and identify robust topological features---e.g., 
those that correspond to long-lived bars in a barcode
diagram.  Since edges of the witness complex are both created and
destroyed with increasing reconstruction dimension, the idea of
persistence would not seem to work for this case (except perhaps in
the sense of zigzag persistence \cite{Carlsson10a}).  However, if we
use $\Delta m$-persistent edges---edges that exist for some
\emph{range} of $m$---we can get a filtration.  If $L^s(\Delta m)$ is the
set of edges that exist for a range $\Delta m$, then $L^s(\Delta m)
\subset L^s(\Delta m-1)$.  The same would be true for a (clique) complex
built from these edges.  Note that this gives, somewhat surprisingly,
an inclusion as $\Delta m$ decreases, the reverse of what one might
think at first.

%%%%%%%%%%%%%%%%%%%%
\section{Conclusion \& Future Work}
\label{sec:concl}

We have shown that it is possible compute the topology of an invariant set
of a dynamical system using a coarse-grained simplicial complex---a
\emph{witness complex}---built from a low-dimensional reconstruction
of a scalar time series.  These results have a
number of interesting implications.  Among other things, they suggest
that the traditional delay-coordinate reconstruction process is
excessive if one is only interested in topological structure.  
Indeed, this explains why it is possible to construct
accurate predictions of the future state of a high-dimensional
dynamical system using a two-dimensional delay-coordinate
reconstruction \cite{josh-chaos15}.  The delay-coordinate machinery
strives to obtain a diffeomorphism---not a homeomorphism---between the true and
reconstructed attractors.  However, many of the important properties of
attractors (transitivity, continuity, recurrence, entropy, etc.) are
topological, so requiring only a homeomorphism is natural
and more efficient \cite{mischaikow99}.

Given a finite set of data, and finite computational power, one can never
compute the \emph{full} topology, of course. Nevertheless, it is useful to
obtain a coarse approximation of the topological features---e.g., the two
main holes in the canonical Lorenz attractor.  Moreover, the
scale of the exploration is controlled by the spacing between the landmarks in
the witness complex, so finer structure can be observed if needed.  There is a
computational cost, of course,
even with the natural parsimony of the witness complex.  

To choose the various free parameters in this approach, we used ideas from
\emph{persistent homology} to investigate
the dependence of the complex on the number
of landmarks $\ell$ and on the fuzziness parameter $\eps$.  It might be possible
to make this multi-parameter
exploration of persistence rigorous using ideas related to those
of \cite{Carlsson10b}.  To study persistence across reconstruction
dimension we introduced the maximum lifespan $\Delta m$. 
This parameter gives a filtration so that the standard persistent homology
methods apply. We found that the long-lived edges determine the 
core of the complex and that the shorter-lived edges fill in the
fine-grained structure. Of course, some edges are superfluous for the 
homology, and destroying them does not change the Betti numbers.
This is part of the strength of the witness complex
approach, as well as the noise immunity that we have observed.

In this paper we have used a fuzzy witness relation based on \Eq{witnessSet} 
and built a clique complex \Eq{fuzzyComplex}.
Since the structure of such a complex is determined by its edges, it is computationally 
efficient. In the future we plan to investigate other choices for the witness
relation, e.g. \cite{deSilva04,Attali07}, and to see if removing the clique 
assumption allows for greater fidelity without excess computational burden. It would be also be interesting to investigate how the robustness with respect to noise varies with the 
choice of relation.

Our approach not only provides a strategy for selecting a reconstruction
dimension that reveals the (approximate) homology of the attractor; it is also a
step along the path to detecting and characterizing bifurcations \cite{BerwaldGV13}.
Suppose that, for example, the data corresponds trivially to an
equilibrium, that is to a set with $\beta_k = 0$ for all $k>0$, but
that this equilibrium undergoes a bifurcation to an oscillatory
regime. In this case, a shift in $\beta_1$ signals a regime change.
Regime changes in a nonstationary data set could also be due dynamical parameters
changing with time, or to switching from one dynamical system to another, e.g., such
as an iterated function system \cite{Alexander12}.

Further analysis of the dynamics, beyond
the static homology of the points in the reconstructed trajectory,
requires the construction of a map on the witness complex that is
induced by the temporal shift on the time series.
Following the ideas of \cite{Alexander15}, the time-ordering of the data gives rise to a 
simplicial multi-valued map, that we call the \emph{witness map},
and under certain conditions this map
can be shown to induce a map on homology that is the same as any
continuous ``selector" and allows one to compute the Conley index
which can be used to prove the existence of various invariant sets.

This paper uses a single example---the classic Lorenz system---but we
have observed similar results in other dynamical systems. We plan in the 
future to do a careful exploration of additional systems, both maps and
flows, to explore the effects of the size and spacing of the
data, to study the interaction with the fuzziness parameter $\eps$
and the number $\ell$ of landmarks.  The edge-length distribution might
be a useful aid in choosing good values for the former; for the
latter, it may be useful to examine the distribution of
witness-landmark distances.  The underlying problem is not simple; the
data sample the invariant set, the landmarks sample the data, and the complex
reflects the geometric relationships between the landmarks and the
rest of the data.  One could get at some of this by changing the data
length---say, using half of the time-series data---repeating the
analysis, and looking to see if the ``best'' parameter values change.
It will be important to determine guidelines for the optimal number and 
distribution of the landmarks, and how this
should scale with $m$ for a fixed trajectory.
%or whether distributing landmarks evenly in time is the best way to go.  

\section*{Acknowledgments}

This material is based upon work supported by the National Science
Foundation under Grants \#CMMI-1245947, \#CNS-0720692, and
\#DMS-1211350. Any opinions, findings, and conclusions or
recommendations expressed in this material are those of the authors
and do not necessarily reflect the views of the National Science
Foundation. The authors thank the IMA for hosting the Algebraic
Topology and Dynamics conference in February 2014, and K. Mischaikow
for numerous useful discussions.

\bibliographystyle{elsarticle-num}
\bibliography{master-refs}

\end{document}